\documentclass[11pt]{article}
\usepackage{latexsym}
\usepackage{amssymb}

\def\qed{\hfill$\Box$}

\newcommand{\bydef}{\mbox{$ \;\stackrel{\triangle}{=}\; $}}

\newcommand{\leqa}{\mbox{$ \;\stackrel{(a)}{\leq}\; $}}
\newcommand{\leqb}{\mbox{$ \;\stackrel{(b)}{\leq}\; $}}

\newcommand{\geqa}{\mbox{$ \;\stackrel{(a)}{\geq}\; $}}
\newcommand{\geqb}{\mbox{$ \;\stackrel{(b)}{\geq}\; $}}
\newcommand{\geqc}{\mbox{$ \;\stackrel{(c)}{\geq}\; $}}

\newcommand{\eqa}{\mbox{$ \;\stackrel{(a)}{=}\; $}}
\newcommand{\eqb}{\mbox{$ \;\stackrel{(b)}{=}\; $}}

\newcommand{\subD}{_{{}_D}}
\newcommand{\subDi}{_{{}_{D_i}}}
\newcommand{\RL}{{\mathbb R}}
\newcommand{\NN}{{\mathbb N}}
\newcommand{\RN}{{\mathbb Q}}
\newcommand{\IND}{{\mathbb I}}
\newcommand{\BBP}{{\mathbb P}}

\newcommand{\PR}{\mbox{\rm Pr}}

\newcommand{\Ahat}{\mbox{$\hat{A}$}}
\newcommand{\Ahatn}{\mbox{$\hat{A}^n$}}
\newcommand{\Ahatk}{\mbox{$\hat{A}^k$}}

\newcommand{\Dmin}{\mbox{$D_{\rm min}$}}
\newcommand{\rhomin}{\mbox{$\rho_{\rm min}$}}
\newcommand{\rhomax}{\mbox{$\rho_{\rm max}$}}
\newcommand{\Lmax}{\mbox{$L_{\rm max}$}}
\newcommand{\Dmax}{\mbox{$D_{\rm max}$}}
\newcommand{\dmax}{\mbox{$d_{\rm max}$}}

\newcommand{\Dminn}{\mbox{$D_{\rm min}^{(n)}$}}
\newcommand{\Dmink}{\mbox{$D_{\rm min}^{(k)}$}}
\newcommand{\Dminone}{\mbox{$D_{\rm min}^{(1)}$}}

\newcommand{\Rmin}{\mbox{$R_{\rm min}$}}
\newcommand{\LA}{\mbox{$\Lambda$}}

\newcommand{\la}{\mbox{$\lambda$}}

\newcommand{\las}{\mbox{\scriptsize$\lambda$}}

\def\be{\begin{eqnarray}}
\def\ee{\end{eqnarray}}
\def\ben{\begin{eqnarray*}}
\def\een{\end{eqnarray*}}

\input{epsf}

\title{Efficient Sphere-Covering and
Converse Measure Concentration\\
Via Generalized Source Coding Theorems}

\author{I. Kontoyiannis}

\date{\today}

\begin{document}
\bibliographystyle{plain}
\maketitle

\thispagestyle{empty}
\setcounter{page}{0}

\vspace{0.3in}

\centerline{\bf Abstract}
{\small
\begin{quote}
Suppose $A$ is a finite set, let $P$ be a 
discrete probability distribution on $A$, 
and let $M$ be an arbitrary ``mass'' 
function on $A$. We give a precise 
characterization of the 
most efficient way in which $A^n$ can be 
almost-covered using spheres of a fixed 
radius. An almost-covering is a subset 
$C_n$ of $A^n$, such that the union of 
the spheres centered at the points of 
$C_n$ has probability close to one with 
respect to the product distribution $P^n$. 
Spheres are defined in terms of a 
single-letter distortion measure on $A^n$,
an efficient covering is one with small 
mass~$M^n(C_n)$, and $n$ is typically large. 
In information-theoretic terms, the sets
$C_n$ are rate-distortion codebooks, but
instead of minimizing their size we seek
to minimize their mass. With different 
choices for $M$ and the distortion measure 
on $A$ our results give various corollaries 
as special cases, including Shannon's 
classical rate-distortion theorem, a version 
of Stein's lemma (in hypothesis testing), 
and a new converse to some measure-concentration 
inequalities on discrete 
spaces. Under mild conditions, we generalize 
our results to abstract spaces and non-product 
measures.

\vspace{0.4in}


{\bf Index Terms} -- Sphere covering, 
measure-concentration,
data compression, large deviations.

\end{quote}
}

\footnotetext[1]{Research 
supported in part by a grant from the Purdue 
Research Foundation.}
\footnotetext[2]{The author is
with the Department
of Statistics, Purdue University,
1399 Mathematical Sciences Building,
W.~Lafayette, IN 47907-1399.
Email: {\tt yiannis@stat.purdue.edu}\,
Web: {\tt www.stat.purdue.edu/$\!\mathtt{\sim}$yiannis}\,
Phone: int+1 (765) 494-6033\,
Fax:   int+1 (765) 494-0558  }

\newpage 

\section{Introduction}

Suppose $A$ is a finite set and let $P$ a
discrete probability mass function on $A$ 
(more general probability spaces are considered 
later).  Assume that the distortion 
(or distance) $\rho(x,y)$ between 
two symbols (or points)
$x$ and $y$ from $A$ 
is measured by a fixed 
$\rho:A\!\times\!A\to[0,\infty)$,
and for each $n\geq 1$ define a 
single-letter distortion measure
(or coordinate-wise distance function)
$\rho_n$ by
\be
\rho_n(x_1^n,y_1^n)=\frac{1}{n}\sum_{i=1}^n\rho(x_i,y_i),
\label{eq:sldistance}
\ee
for $x_1^n=(x_1,x_2,\ldots,x_n)$ and 
$y_1^n=(y_1,y_2,\ldots,y_n)$ in $A^n$.

Given a $D\geq 0$, we want 
to ``almost'' cover the product space
$A^n$ using
a finite number of balls (or ``spheres'') $B(y_1^n,D)$,
where
\be
B(y_1^n,D)=\{x_1^n\in A^n\;:\;\rho_n(x_1^n,y_1^n)\leq D\}
\label{eq:balldef}
\ee
is the (closed) ball of 
distortion-radius $D$ centered 
at $y_1^n\in A^n$. For our purposes,
an ``almost covering'' is a subset $C\subset A^n$, 
such that the union of the balls of radius $D$ centered
at the points of $C$ have large $P^n$-probability,
that is,
\be
P^n\left([C]\subD\right)
\;\;\;\mbox{is close to 1,}
\label{eq:cover}
\ee
where $[C]\subD$ is the $D$-{\em blowup of} $C$
\ben
[C]\subD=
\bigcup_{y_1^n\in C} B(y_1^n,D)
=\{x_1^n\;:\;\rho_n(x_1^n,y_1^n)\leq D\;\;
	\mbox{for some}\;y_1^n\in C\}.
\een
More specifically, given a ``mass function''
$M:A\to(0,\infty),$ we are interested 
in covering $A^n$ {\em efficiently}, namely, 
finding sets $C$ that satisfy (\ref{eq:cover})
and also have small mass
\ben
M^n(C)=\sum_{y_1^n\in C} M^n(y_1^n)
=\sum_{y_1^n\in C}\prod_{i=1}^nM(y_i).
\een
One way to state our main question of 
interest is as follows:
\ben
(*)\;
  \cases{
	\hspace{0.05in}
        \mbox{\em If the sets $\{C_n\;;\;n\geq 1\}$ 
		  asymptotically $D$-cover 
		  $A^n$, that is,}&\cr
	\hspace{1.5in}
		\left.P^n\left([C_n]\subD\right)\to 1\right.
		\;\;\;\mbox{as $n\to\infty$},&\cr
	\hspace{0.05in}
	\mbox{\em how small can their masses $M^n(C_n)$ 
		be?}\cr
	}
	\hspace{0.3in}
\een


\noindent
Question $(*)$ is partly motivated by the
fact that several interesting questions can
be easily restated in this form.
Three such examples are presented below,
and in the remainder of the paper $(*)$ 
is addressed and answered in detail. 
In particular, it is shown that $M^n(C_n)$ 
typically grows (or decays) exponentially 
in $n$, and an explicit lower bound, valid 
for all finite $n$, is given for the 
exponent $(1/n)\log M^n(C_n)$ 
of the mass of an arbitrary $C_n$.
[Throughout the paper, `log' denotes 
the natural logarithm.]
Moreover, a sequence of sets $C_n$ 
asymptotically achieving this lower 
bound is exhibited, showing that 
it is best possible. The outline of
the proofs follows, to some extent,
along similar lines as the proof of 
Shannon's rate-distortion theorem
\cite{shannon:59}.
In particular, the ``extremal''
sets $C_n$ achieving the lower
bound are constructed probabilistically;
each $C_n$ consists of a collection of 
points $y_1^n$ generated by taking 
independent and identically distributed 
(IID) samples from a suitable
distribution on $A^n$, but (unlike Shannon)
here we need to condition on seeing 
typical realizations, making the 
individual elements of the random 
$C_n$ non-IID.

\vspace{0.1in}

{\sc Example 1. (Measure 
Concentration on the Binary Cube)}
Take $A=\{0,1\}$ so that $A^n$ is the
$n$-dimensional binary cube consisting
of all binary strings of length $n$,
and let $P^n$ be a product probability
distribution on $A^n$.
Write $\rho_n(x_1^n,y_1^n)$ for the 
normalized Hamming distortion
between $x_1^n$ and $y_1^n$,
so that $\rho_n(x_1^n,y_1^n)$
is the proportion of
mismatches between the two strings;
formally:
\be
\rho_n(x_1^n,y_1^n)=\frac{1}{n}\sum_{i=1}^n
\IND_{
\{x_i\neq y_i\}
},\;\;\;\;
x_1^n, y_1^n\in A^n.
\label{eq:hammingdist}
\ee
Geometrically, if $A^n$ is 
given the usual nearest-neighbor
graph structure (two points are connected if and 
only if they differ in exactly one coordinate),
then $\rho_n(x_1^n,y_1^n)$ is the
graph distance between $x_1^n$ and 
$y_1^n$, normalized by $n$.

A well-known measure-concentration inequality
for subsets $C_n$ of $A^n$ states that,
for any
$D\geq 0$,
\be
P^n([C_n]\subD)\geq 1 -\frac{e^{-nD^2/2}}{P^n(C_n)}.
\label{eq:mcdiarmid}
\ee
[See
Proposition~2.1.1 in the
comprehensive account by
Talagrand \cite{talagrand:95},
or Theorem~3.5 in 
the review paper by
McDiarmid \cite{mcdiarmid:98},
and the references therein.] Roughly
speaking,
(\ref{eq:mcdiarmid}) says that ``if
$C_n$ is not too small, $[C_n]\subD$ is almost 
everything.'' In particular, it implies
that for any sequence of
sets $C_n\subset A^n$
and any $D\geq 0$,
\be
\mbox{\em if}\;\;\;\;\;
\liminf_{n\to\infty}\frac{1}{n}\log P^n(C_n) > -D^2/2,
\;\;\;\;\;\mbox{\em then}\;\;\;\;\;
P^n([C_n]\subD)\to 1.
\label{eq:directMC}
\ee
A natural question to ask is
whether there is a converse
to the above statement: If 
$P^n([C_n]\subD)\to 1$, how small
can the probabilities of the $C_n$ be?
Taking $M\equiv P$, this reduces to 
question $(*)$ above. 
In this context, $(*)$
can be thought of as the opposite
of the usual isoperimetric problem.
We are looking for sets with 
the ``largest possible boundary'';
sets $C_n$ whose $D$-blowups (asymptotically)
cover the entire space, but whose volumes
$P^n(C_n)$ are as small as possible.
A precise answer for this problem is given 
in Corollary~3 and the discussion following 
it, in the next section.

\vspace{0.1in}

{\sc Example 2. (Lossy Data Compression)}
Let $A$ be a finite alphabet so that 
$A^n$ consists of all possible messages
of length $n$ from $A$, and assume that 
messages are generated by a memoryless
source, with distribution $P^n$ on $A^n$. 
A code for these messages consists of 
a codebook $C_n\subset A^n$
and an encoder $\phi_n:A^n\to C_n$.
If we think of $\rho_n(x_1^n,y_1^n)$ as
the distortion between a message
$x_1^n$ and its reproduction $y_1^n$,
then for any given codebook $C_n$ the 
best choice for the encoder is clearly 
the map $\phi_n$ taking each $x_1^n$ 
to the $y_1^n$ in $C_n$ which minimizes
the distortion $\rho_n(x_1^n,y_1^n)$.
Hence, at least conceptually,
finding good codes is the
same as finding good codebooks.
More specifically, if $D\geq 0$ is the 
maximum amount of distortion we are 
willing to tolerate, then a sequence
of good codebooks $\{C_n\}$ is one 
with the following properties:
\begin{itemize}
\item[$(a)$] The probability
	of encoding a message with distortion
	exceeding $D$ is asymptotically negligible:
	$$P^n([C_n]\subD)\to 1.$$
\item[$(b)$] Good compression is achieved,
	that is, the sizes $|C_n|$ of the 
	codebooks are small.
\end{itemize}

\noindent
What is the best achievable compression
performance? That is, if the
codebooks $\{C_n\}$ satisfy $(a)$, how 
small can their sizes be? Shannon's 
classical source coding theorem (cf. 
\cite{shannon:59}\cite{berger:book}) 
answers this question. In our notation, 
taking $M\equiv 1$ reduces the question 
to a special case of $(*)$, and in 
Corollary~2 in the next section
we recover Shannon's 
theorem as a special case of 
Theorems~1 and~2.


\vspace{0.1in}

{\sc Example 3. (Hypothesis Testing)}
Let $A$ be a finite set and $P_1,\,P_2$ be 
two probability distributions on $A$.
Suppose that the null hypothesis that a sample
$X_1^n=(X_1,X_2,\ldots,X_n)$ of $n$
independent observations comes from $P_1$ 
is to be tested against the simple 
alternative hypothesis that $X_1^n$ comes 
from $P_2$. A test between these two
hypotheses can be thought of as a 
decision region $C_n\subset A^n$:
If $X_1^n\in C_n$ we declare that 
$X_1^n\sim P_1^n$, otherwise
we declare $X_1^n\sim P_2^n$.
The two probabilities of error associated with 
this test are
\be
\alpha_n = P_1^n(C_n^c)
\;\;\;\;\;
\mbox{and}
\;\;\;\;\;
\beta_n  = P_2^n(C_n).
\label{eq:errorprob}
\ee
A good test has these two probabilities
vanishing as fast as possible, and
we may ask, if $\alpha_n\to 0$,
how fast can $\beta_n$ decay to zero?
Taking $\rho$ to be Hamming distortion, $D=0$,
$P=P_1$, and $M=P_2$, this reduces to our original 
question $(*)$. In Corollary~1 in the next section
we answer this question by deducing
a version of Stein's lemma from 
Theorems~1 and~2. It is worth noting 
that the connection between
questions in hypothesis testing and information 
theory goes at least as far back as Strassen's 
1964 paper \cite{strassen:64b} (see also Blahut's 
paper \cite{blahut:74} in 1974, and Csisz{\'{a}}r 
and K{\"{o}}rner's book \cite{csiszar:book} for a 
detailed discussion). 

\vspace{0.1in}

The rest of the paper is organized as follows.
In Section~2, Theorems~1 and~2 provide an 
answer to question $(*)$. 
In the remarks and corollaries following 
Theorem~2 we discuss and interpret this answer, 
and we present various applications along the lines
of the three examples above. Theorem~1 is proved
in Section~2 and Theorem~2 is proved in Section~3.
In Section~4 we consider the same problem in a 
much more general setting. We let $A$ be an 
abstract space, and instead of product 
measures $P^n$ we consider the $n$-dimensional 
marginals $P_n$ of a stationary measure 
$\BBP$ on $A^\NN$. 
In Theorems~3 and~4 we give analogs 
of Theorems~1 and~2, which hold essentially
as long as the spaces $(A^n,P_n)$ can be 
almost-covered by countably many 
$\rho_n$-balls.
Although the results of Section~2 are 
essentially subsumed by Theorems~3 and~4,
it is possible to give simple, elementary 
proofs for the special case treated in 
Theorems~1 and~2, so we give separate 
proofs for these results first. The more 
general Theorems~3 and~4 are proved in 
Section~5, and the Appendix contains the 
proofs of various technical steps needed 
along the way.

\section{The Discrete Memoryless Case}

Let $A$ be a finite set and $P$ be a 
discrete probability mass function on $A$. 
Fix a $\rho:A\!\times\!A\to[0,\infty)$,
and for each $n\geq 1$ let 
$\rho_n$ be the corresponding
single-letter distortion measure
(or coordinate-wise distance function)
on $A^n$ defined as in (\ref{eq:sldistance}).
Also let $M:A\to(0,\infty)$ be an 
arbitrary positive mass function on $A$.
We assume, without loss of generality, 
that $P(a)>0\,$ for all $a\in A$, and
also that for each $\,a\in A$ there exists
a $\,b\in A$ with $\rho(a,b)=0$ (otherwise
we may consider 
$\rho'(x,y)=[\rho(x,y)-\min_{z\in A}\rho(x,z)]$
instead of $\rho$). Let $\{X_n\}$ denote a sequence of IID
random variables with distribution $P$, and
write $\BBP=P^\NN$ for the
product measure on $A^\NN$ equipped with 
the usual $\sigma$-algebra generated by 
finite-dimensional cylinders. We write
$X_i^j$ for vectors of random variables
$(X_i,X_{i+1},\ldots,X_j)$, 
$1\leq i\leq j\leq\infty$, and similarly 
$x_i^j=(x_i,x_{i+1},\ldots,x_j)\in A^{j-i+1}$
for realizations of these random variables.

Next we define the rate function $R(D)$ that 
will provide the lower bound on 
the exponent of the mass of an 
arbitrary $C_n\subset A^n$. For $D\geq 0$ and $Q$ 
a probability measure on $A$, let
\be
I(P,Q,D)=\inf_{W\in{\cal M}(P,Q,D)}
	H(W\|P\!\times\!Q)
\label{eq:Idef}
\ee
where $H(\mu\|\nu)$ denotes the relative entropy
between two discrete probability mass functions 
$\mu$ and $\nu$ on a finite set $S$,
\ben
H(\mu\|\nu)=\sum_{s\in S}\mu(s)
\log\frac{\mu(s)}{\nu(s)},
\een
and where ${\cal M}(P,Q,D)$
consists of all probability measures
$W$ on $A\!\times\!A$ 
such that $W_X$, the first 
marginal of $W$, is equal to $P$, $W_Y$, the
second marginal, is $Q$, and $E_W[\rho(X,Y)]\leq D$;
if ${\cal M}(P,Q,D)$ is empty,
we let $I(P,Q,D)=\infty$.
The rate function $R(D)$ is defined by
\be
R(D)=R(D;P,M) = \inf_{Q}
		\left\{ I(P,Q,D) + E_Q[\log M(Y)]
		\right\},
\label{eq:simpleRdef}
\ee
where the infimum is over all 
probability distributions $Q$ on $A$.
Recalling the definition of the
{\em mutual information} between two
random variables, $R(D)$ can equivalently
be written in a more information-theoretic
way. If $(X,Y)$ are random variables
(or random vectors) with joint distribution 
$W$ and corresponding marginals $W_X$ and 
$W_Y$, then the mutual information between 
$X$ and $Y$ is defined as
$$I(X;Y)=H(W\|W_X\!\times\!W_Y).$$
Combining the two infima in (\ref{eq:Idef})
and (\ref{eq:simpleRdef}) we can write
\be
R(D)=
\inf_{(X,Y):\;X\sim P,\;E\rho(X,Y)\leq D} \left\{
	I(X;Y) + E[\log M(Y)] \right\}
\label{eq:infoRdef}
\ee
where the infimum is taken over
all jointly distributed 
random variables $(X,Y)$
such that $X$ has distribution $P$
and $E\rho(X,Y)\leq D$.
For any $x_1^n\in A^n$ and $C_n\subset A^n$, write
\ben
\rho_n(x_1^n,C_n)=\min_{y_1^n\in C_n}\rho_n(x_1^n,y_1^n).
\een

In the following two Theorems we answer
question $(*)$ stated in the Introduction. 
Theorem~1 contains a lower bound
(valid for all finite $n$)
on the mass of an arbitrary 
$C_n\subset A^n$,
and Theorem~2 shows that this 
bound is asymptotically tight.
In information-theoretic terms,
Theorems~1 and~2 can be thought of 
as generalized direct and converse 
coding theorems, for minimal-mass
(rather than minimal-size) codebooks.

\vspace{0.1in}

{\sc Theorem 1.} {\em Let $C_n\subset A^n$ 
be arbitrary and write $D=E_{P^n}[\rho_n(X_1^n, C_n)]$.
Then
\ben
\frac{1}{n}\log M^n(C_n)\geq R(D).
\een
}

\vspace*{-0.15in}

{\sc Theorem 2.} {\em Assume that $\rho(x,y)=0$ if
and only if $x=y$. 
For any $D\geq 0$ 
and any $\epsilon>0$ 
there is a sequence of sets $\{C_n\}$
such that:
\ben
&(i)&\hspace{0.2in}
	\frac{1}{n}\log M^n(C_n)\leq R(D)+\epsilon
	\hspace{0.3in} \mbox{for all $n\geq 1$}\\
&(ii)&\hspace{0.2in}
	\rho_n(X_1^n,C_n)\leq D
	\hspace{0.3in} 
	\mbox{eventually, }
	\BBP-a.s.
\een
}

\vspace*{-0.1in}

{\sc Remark 1.}
Part $(ii)$ of Theorem~2 
says that $\IND_{[C_n]\subD}(X_1^n)\to 1$
with probability one, so
by Fatou's lemma, 
$P^n\left([C_n]\subD\right)\to 1.$
From this and $(i)$ it is easy to deduce
the following
alternative version 
of Theorem~2 (see the Appendix for a proof):
{\em
For any $D\geq 0$ 
there is a sequence of sets $\{C^*_n\}$
such that:}
\ben
&(i')&\hspace{0.2in}
	\limsup_{n\to\infty}\;
	\frac{1}{n}\log M^n(C^*_n)\leq R(D)\\
&(ii')&\hspace{0.2in}
	P^n([C^*_n]\subD)\to 1,\;\;\;\mbox{and}\\
&(iii')&\hspace{0.2in}
	\limsup_{n\to\infty}\;
	E_{P^n}[\rho_n(X_1^n,C^*_n)]\leq D
\een

{\sc Remark 2.} As will become evident from the proof
of Theorem~2, the additional assumption on $\rho$ is only
made for the sake of simplicity, and it is not necessary
for the validity of the result. In particular, it allows
us to give a unified argument for the cases $D=0$ and $D>0$.

\vspace{0.1in}

Theorem~1 is proved at the end of this section,
and Theorem~2 is proved in Section~3. Although
the proof of Theorem~2 is somewhat technical,
the idea behind the construction of the extremal 
sets $C_n$ is simple: Suppose $Q^*$ is a
probability measure on $A$ achieving the infimum 
in the definition of $R(D)$, so that
	$$R(D)=I(P,Q^*,D)+E_{Q^*}[\log M(Y)]\bydef I^*+L^*.$$ 
Write $Q_n^*$ for the product measure $(Q^*)^n$, 
and let $\widehat{Q}_n$ be the measure obtained by 
conditioning $Q_n^*$ to the set of points 
$y_1^n\in A^n$ whose empirical measures 
(``types'') are uniformly close to $Q^*$. 
Then let $C_n$ consist of approximately
$e^{nI^*}$ points $y_1^n$ drawn
IID from $\widehat{Q}_n$. Each point in the
support of $\widehat{Q}_n$ has mass 
$M^n(y_1^n)\approx e^{nL^*}$ and $C_n$
contains about $e^{nI^*}$ of them, so 
$M^n(C_n)$ is close to 
$e^{nI^*}e^{nL^*}=e^{nR(D)}.$ The main
technical content of the proof is therefore
to prove $(ii)$, namely, that $e^{nI^*}$
points indeed suffice to almost 
$D$-cover $A^n$.

The above construction also provides a nice
interpretation for $R(D)$. If we had 
started with a different measure $Q$ 
in place of $Q^*$, we would have ended
up with sets $C'_n$ of size 
$\mbox{$\approx \exp(nI(P,Q,D))$}$, consisting of
points $y_1^n$ of mass 
$M^n(y_1^n)\approx \exp(nE_{Q}(\log M(Y))),$
and the total mass of $C'_n$ would be
$$M^n(C'_n)\approx
	\exp\{n[I(P,Q,D)+E_{Q}(\log M(Y))]\}.$$
By optimizing over the choice of $Q$ in 
(\ref{eq:simpleRdef}) we are balancing 
the tradeoff between the size and the 
weight of the set $C_n$, between a few 
heavy points and many light ones.

It is also worth noting that the extremal
sets $C_n$ above were constructed by taking
samples $y_1^n$ from the {\em non}-product 
measure $\widehat{Q}_n$. Unlike in Shannon's 
proof of the data compression theorem,
here we cannot get away by simply using
the product measure $Q_n^*$. This is 
because we are not just interested 
in how many points $y_1^n$ are needed 
to almost cover $A^n$, but also we need
control their masses $M^n(y_1^n)$. 
Since exponentially many $y_1^n$'s are
required to cover $A^n$, if they are
generated from $Q_n^*$ then there are
bound to be some atypically heavy ones,
and this drastically increases the total mass 
$M^n(C_n)$. Therefore, by restricting $Q_n^*$ 
to be supported on the set of $y_1^n\in A^n$ 
whose empirical measures are uniformly close 
to $Q^*$, we are ensuring that the masses 
of the $y_1^n$ will be essentially constant,
and all approximately equal to $e^{nL^*}$.

Next we derive corollaries from 
Theorems~1 and~2, along the lines 
of the examples in the Introduction.  
First, in the context of hypothesis 
testing, let $P_1,\,P_2$ be two 
probability distributions on $A$
with all positive probabilities.
Suppose that the null hypothesis that
$X_1^n\sim P_1^n$ is to be tested against the
alternative $X_1^n\sim P_2^n$. Given a test
with an associated decision region 
$C_n\subset A^n$, its two probabilities 
of error $\alpha_n$ and $\beta_n$ are 
defined as
in (\ref{eq:errorprob}). In the notation
of this section, let $\rho_n$ be Hamming 
distortion as in (\ref{eq:hammingdist}), 
$P=P_1$ and $M=P_2$. Observe that, here,
$$E_{P_1^n}[\rho_n(X_1^n,C_n)]
	\leq E_{P_1^n}[\IND_{C_n^c}(X_1^n)]
		= P_1^n(C_n^c),$$
and define, in the notation of
(\ref{eq:simpleRdef}), the error
exponent
$$\varepsilon(\alpha)=-R(\alpha;P_1,P_2),
\;\;\;\;\alpha\in[0,1].$$
Noting that $\varepsilon(0)=H(P_1\|P_2)$,
from Theorems~1 and~2 and Remark~1
we obtain the following
version of Stein's lemma (see 
Lemma~6.1 in Bahadur's monograph 
\cite{bahadur:book}, or Theorem~12.8.1
in \cite{cover:book}).

\vspace{0.1in}

{\sc Corollary 1. (Hypothesis Testing)}
{\em 
Let $\alpha=\alpha_n=P_1^n(C_n^c)$ and
$\beta=\beta_n=P_2^n(C_n)$ be the two
types of error probabilities associated
with an arbitrary sequence of tests $\{C_n\}$.
\begin{itemize}
\item[$(i)$] For all $n\geq 1$,
$\,\beta\geq e^{-n\varepsilon(\alpha)}.$
\item[$(ii)$] If $\alpha_n\to 0$, then
\ben	
\liminf_{n\to\infty}\frac{1}{n}\log\beta_n\geq -H(P_1\|P_2).
\een
\item[$(iii)$]
There exists a sequence of 
decision regions $C_n$ with associated
tests whose error probabilities
achieve $\alpha_n\to 0$ and
$(1/n)\log\beta_n\to -H(P_1\|P_2),$
as $n\to\infty$.
\end{itemize}
}

Note that, although the
decision regions $C_n$ in $(iii)$ 
above achieve the best exponent in 
the error probability, they are not 
the overall optimal decision regions
in the Neyman-Pearson sense.

In the case of data compression, 
we have random data $X_1^n$
generated by some product distribution $P^n$. 
Given a single-letter distortion measure
$\rho_n$ and a 
maximum allowable distortion level 
$D\geq 0$, our objective is to find
good codebooks $C_n$.
As discussed in Example~2 above,
good codebooks are those 
that asymptotically cover $A^n$,
i.e., $P^n([C_n]\subD)\to 1,$
and whose sizes $|C_n|$ are relatively
small. In our notation, if we take
$M(\cdot)\equiv 1$, then 
$M^n(C_n)=|C_n|$ and 
the rate function $R(D)$ 
(from (\ref{eq:simpleRdef}) 
or (\ref{eq:infoRdef}))
reduces to Shannon's 
{\em rate-distortion function}
\ben
R_S(D)&=&\inf_Q\inf_{W\in{\cal M}(P,Q,D)} H(W\|P\!\times\!Q)\\
&=& \inf_{(X,Y):\;X\sim P,\;E\rho(X,Y)\leq D} I(X;Y).
\een
From Theorems~1 and~2 and Remark~1
we recover Shannon's source coding theorem
(see \cite{shannon:59}\cite{berger:book}).

\vspace{0.1in}

{\sc Corollary 2. (Data Compression)}
{\em For any $n\geq 1$, if the average 
distortion achieved by a codebook $C_n$ is 
$D=E_{P^n}[\rho_n(X_1^n,C_n)]$, then
$$\frac{1}{n}\log |C_n|\geq R_S(D).$$
Moreover, for any $D\geq 0$, there is a 
sequence of codebooks $\{C_n\}$ such that
$E_{P^n}[\rho_n(X_1^n,C_n)]\to D,$ the 
codebooks $C_n$ asymptotically cover $A^n$,
$ P^n([C_n]\subD)\to 1$, and
$$\lim_{n\to\infty}\;
        \frac{1}{n}\log |C_n| = R_S(D).$$
}

%

Finally, in the context of measure-concentration, 
taking $M=P$ 
and writing $R_C(D)$ for the 
concentration exponent $R(D;P,P)$, 
we get:

\vspace{0.1in}

{\sc Corollary 3. (Converse Measure Concentration)}
{\em Let $\{C_n\}$ be arbitrary sets.
\begin{itemize}
\item[$(i)$] For any $n\geq 1$, if $D=E_{P^n}[\rho_n(X_1^n,C_n)]$,
then $\,P^n(C_n)\geq e^{nR_{C}(D)}.$
\item[$(ii)$]
If $\,P^n([C_n]\subD)\to 1,$ then 
$$\liminf_{n\to\infty}\frac{1}{n}
\log P^n(C_n)\geq R_{C}(D).$$
\item[$(iii)$]
There is a sequence of sets
$\{C_n\}$ such that 
$\,P^n([C_n]\subD)\to 1$
and
$(1/n) \log P^n(C_n)\to R_{C}(D),$
as $n\to\infty$.
\end{itemize}
}

In particular, in the case of the 
binary cube, part $(ii)$ of the corollary 
provides a precise converse to the 
measure-concentration statement in 
(\ref{eq:directMC}). Although the 
concentration exponent $R_C(D)=R(D;P,P)$
is not as explicit as the exponent
$\,-D^2/2\,$ in (\ref{eq:directMC}), 
$R_C(D)$ is a well-behaved function
and it is easy to evaluate it 
numerically. For example, Figure~1 
shows the graph of
$R_C(D)$ in the case of the binary 
cube, with $P$ being the Bernoulli
measure with $P(1)=0.4$. 
Various easily checked properties 
of $R(D)=R(D;P,M)$ are stated in 
Lemma~1, below; proof outlines 
are given in the Appendix.

\begin{figure}[ht]
\centerline{\epsfxsize 3in \epsfbox{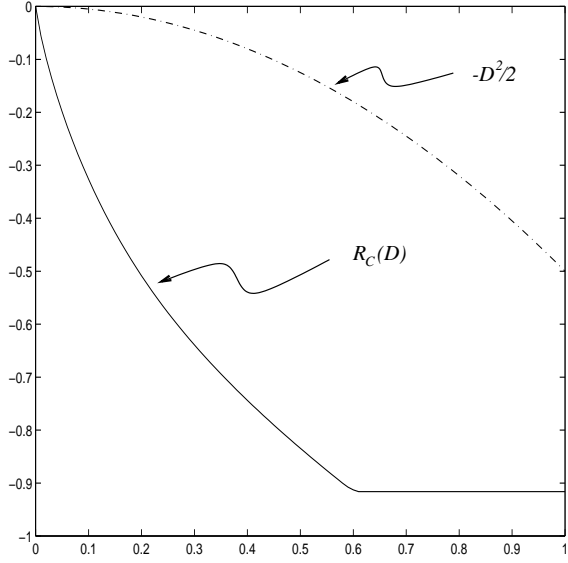}}
\caption{Graph of the function $R_C(D)=R(D;P,P)$ for
$0\leq D\leq 1$, in the case of the binary cube
$A^n=\{0,1\}^n$, with $P(1)=0.4$.}
\end{figure}

As mentioned in the Introduction,
the question considered in Corollary~3 
can be thought of as the opposite 
of the usual isoperimetric problem.
Instead of large sets with small 
boundaries, we are looking for 
{\em small} sets with the
{\em largest possible boundary.} It is
therefore not surprising that 
the extremal sets in (\ref{eq:directMC})
and in Corollary~3 are very different.
In the classical isoperimetric
problem, the extremal sets typically 
look like Hamming balls around 
$0^n=(0,0,\ldots,0)\in A^n$,
$B_n=\{x_1^n\;:\;\rho_n(x_1^n,0^n)\leq r/n\}$
(see the discussions in
Section~2.3 of \cite{talagrand:95}, 
p.~174 in \cite{mcdiarmid:89}, 
or the original paper by Harper 
\cite{harper:66}), while the extremal
sets in our case are collections
of vectors $y_1^n$ drawn IID
from the measure $\widehat{Q}_n$ 
on $A^n$.


\newpage

{\sc Lemma 1.}{\em
$(i)$
The infima in the definitions of $R(D)$
and $I(P,Q,D)$ in (\ref{eq:simpleRdef})
and (\ref{eq:Idef})
are in fact minima.

$(ii)$
$R(D)$ is finite for all $D\geq 0$,
it is nonincreasing and convex in $D$,
and therefore also continuous.

$(iii)$
For fixed $P$ and $Q$, 
$I(P,Q,D)$ is nonincreasing and convex
in $D$, and therefore it is
continuous except possibly 
at the point 
$D=\inf\{D\geq 0\;:\;I(P,Q,D)<\infty\}$.

$(iv)$
If the random variables $X_1^n=(X_1,\ldots,X_n)$
are IID, then for any random vector $Y_1^n$
jointly distributed with $X_1^n$:
\ben
I(X_1^n;Y_1^n)\geq \sum_{i=1}^nI(X_i;Y_i).
\een

$(v)$ If we let $\Rmin=\min\{\log M(y)\;:\;y\in A\}$
and
\ben 
\Dmax=\Dmax(P)=
	\min\{E_P[\rho(X,y)]\;:\;y\;\mbox{such that}\;\log M(y)=\Rmin\},
\een
then
\ben
R(D)\;\mbox{is}\;
	\cases{
        =\Rmin& for $D\geq\Dmax$\cr
        >\Rmin& for $0\leq D<\Dmax$.\cr
              }
\een
}

Next we prove Theorem~1. It is perhaps
somewhat surprising that the proof is 
very short and completely elementary,
relying only on 
Jensen's inequality
and the convexity of $R(D)$.

\vspace{0.1in}

{\em Proof of Theorem~1}:
Given an arbitrary $C_n$,
let $\phi_n:A^n\to C_n$
be a function that maps each $x_1^n\in A^n$
to the closest $y_1^n$ in $C_n$, i.e.,
$\rho_n(x_1^n,\phi(x_1^n))=\rho_n(x_1^n,C_n)$.
For $X_1^n\sim P^n$ let
$Y_1^n=\phi_n(X_1^n)$, 
write $Q_n$ for the
distribution of $Y_1^n$,
and
$W_n(x_1^n,y_1^n) = 
P^n(x_1^n)\IND_{\{\phi_n(x_1^n)\}}(y_1^n)$
for the joint distribution
of $(X_1^n,Y_1^n)$.
Then 
\be
E_{W_n}[\rho_n(X_1^n,Y_1^n)]=D
\label{eq:phidist}
\ee
and by Jensen's inequality, 
\ben
\log M^n(C_n)
&=&\log \left[\sum_{y_1^n\in C_n}\left(Q_n(y_1^n)
	\frac{M^n(y_1^n)}{Q_n(y_1^n)}
	\right)\right]\\
&\geq& \sum_{y_1^n\in C_n}Q_n(y_1^n)\log
	\frac{M^n(y_1^n)}{Q_n(y_1^n)}\\
&=& \sum_{x_1^n,y_1^n\in A^n}W_n(x_1^n,y_1^n)\log
	\frac{W_n(x_1^n,y_1^n)}{P^n(x_1^n)Q_n(y_1^n)}
	+\sum_{y_1^n\in C_n}Q_n(y_1^n)\log M^n(y_1^n).
\een
By the definition of mutual information this equals
$$I(X_1^n;Y_1^n)+E_{Q_n}[\log M^n(Y_1^n)],$$
which, by Lemma~1~$(iv)$, is bounded below by
$$\sum_{i=1}^n\left[I(X_i;Y_i)+E_{Q_n}[\log M(Y_i)]\right].$$
Finally, by the definition of $R(D)$ and its convexity
this is bounded below by
\ben
\sum_{i=1}^n R\left(E_{W_n}[\rho(X_i,Y_i)]\right)
\geq nR\left(\frac{1}{n}\sum_{i=1}^nE_{W_n}[\rho(X_i,Y_i)]\right)
=nR(D)
\een
where the last equality follows from (\ref{eq:phidist}).
\qed

\section{Proof of Theorem 2.}

Let $P$, $D\geq 0$ be fixed,
and $\epsilon>0$ be given.
By Lemma~1~$(i)$ we can pick $Q^*$ and $W^*$
in the definition of $R(D)$ and
$I(P,Q^*,D)$, respectively, such that 
\ben
R(D)
=H(W^*\|P\!\times\!Q^*)+E_{Q^*}[\log M(Y)]
\bydef I^*+L^*.
\een
For $n\geq 1$, write $Q_n^*$ for the product measure
$(Q^*)^n$, and for $y_1^n\in A^n$ let
$$\hat{P}_{y_1^n}=\frac{1}{n}\sum_{i=1}^n
	\delta{y_i}$$
denote the empirical measure of $y_1^n$.
Pick $\delta>0$ (to be chosen later) and
define, for each $n\geq 1$, the set of ``good''
strings
$${\cal G}_n=\{y_1^n\in A^n\;:\;
	\hat{P}_{y_1^n}(b)\leq Q^*(b)+\delta,
	\;\;\forall\, b\in A\}$$
(if ${\cal G}_n$ as defined above is empty -- 
this may only happen for finitely many $n$ -- 
simply let ${\cal G}_n$ consist of a single 
vector $(a,a,\ldots,a)$, with $a\in A$
chosen so that $\log M(a)=\Rmin$).
Also, let $\widehat{Q}_n$ be the measure $Q_n^*$
conditioned on ${\cal G}_n$:
$$\widehat{Q}_n(F)=\frac{Q_n^*(F\cap {\cal G}_n)}{Q_n^*({\cal G}_n)};
\;\;\;\;F\subset A^n.$$
For $n\geq 1$, let
$\{Y(i)=(Y_1(i),Y_2(i),\ldots,Y_n(i))\;;\;i\geq 1\}$
be an IID sequence of random vectors 
$Y(i)\sim \widehat{Q}_n$, and define 
$C_n$ as the collection of the first
$e^{n(I^*+\epsilon/2)}$ of them:
$$C_n=\{Y(i)\;:\;1\leq i\leq e^{n(I^*+\epsilon/2)}\}.$$
By the definition of ${\cal G}_n$, any $y_1^n\in {\cal G}_n$ has
\ben
\frac{1}{n}\log M^n(y_1^n)
= \sum_{b\in A} \hat{P}_{y_1^n}(b)\log M(b)
\leq L^* + \delta\left(
	\sum_{b\in A}\log M(b)
	\right)
\leq L^* + \epsilon/2,
\een
by choosing $\delta>0$ appropriately small.
Therefore,
$$M^n(C_n)\leq e^{n(I^*+\epsilon/2)}e^{n(L^*+\epsilon/2)}
	=e^{n(R(D)+\epsilon)}$$
and $(i)$ of the Theorem is satisfied. 
Let $X_1^n$ be IID random variables
with distribution $P$. To
verify $(ii)$ we will show that 
\be
i_n\leq e^{n(I^*+\epsilon/2)}
\;\;\;\;\mbox{eventually, }\BBP\!\times\!\RN-
	\mbox{a.s.}
\label{eq:i-target}
\ee
where $i_n$ is the index of the first $Y(i)$
that matches $X_1^n$ within $\rho_n$-distortion
$D$,
$$i_n=\inf\{i\geq 1\;:\;\rho_n(X_1^n,Y(i))\leq D\},$$
and $\RN=\prod_{n\geq 1}(\widehat{Q}_n)^\NN$.
Recall the notation
$B(x_1^n,D)=\{y_1^n\in A^n\;:\;\rho_n(x_1^n,y_1^n)\leq D\}.$
For (\ref{eq:i-target}) it suffices to prove 
the following two statements
\be
\limsup_{n\to\infty}\frac{1}{n}\log\left[i_n\,
	\widehat{Q}_n(B(X_1^n,D))\right]
	\leq 0\;\;\;\;
	\BBP\!\times\!\RN-
        \mbox{a.s.}
	\label{eq:4a}\\
\liminf_{n\to\infty}\frac{1}{n}\log
	\widehat{Q}_n(B(X_1^n,D))
	\geq -I^*
	\;\;\;\;
        \BBP-
        \mbox{a.s.}
	\label{eq:4b}
\ee
Proving (\ref{eq:4b}) is the main technical part of
the proof and it will be done last. Assuming it holds,
we will first establish (\ref{eq:4a}). 
For $m\geq 1$ let 
$G_m=\{x_1^\infty\in A^\infty\;:\;
	\widehat{Q}_n(B(x_1^n,D))>0
	\;\;\forall\, n\geq m\},$
and note that by (\ref{eq:4b}),
$\BBP\left(\cup_{m\geq 1}G_m\right)=1.$
Pick $m\geq 1$; for any $n\geq m$, and 
any $x_1^{\infty}\in G_m$, conditional on
$X_1^n=x_1^n$, $i_n$ is a Geometric($p_n$)
random variable with $p_n= \widehat{Q}_n(B(x_1^n,D))$.
So for $\epsilon'>0$ arbitrary
$$\PR\left\{\left.\frac{1}{n}\log \left[i_n\,
        \widehat{Q}_n(B(X_1^n,D))\right]
	>\epsilon'\,\right|\,X_1^n=x_1^n\right\}
	\leq (1-p_n)^{\frac{e^{\epsilon'n}}{p_n}-1}$$
and for all $n$ large enough (independent of
$x_1^n$) this is bounded above by
$$\left[(1-p_n)^{1/p_n}
	\right]^{e^{\epsilon'n-1}}
\leq e^{-e^{\epsilon'n-1}},$$
uniformly over $x_1^\infty\in G_m$. Since
the above right-hand side is summable 
over $n$, by the Borel-Cantelli lemma 
and the fact that $\epsilon'>0$ was arbitrary
we get (\ref{eq:4a}) for $\BBP$-almost all 
$x_1^\infty\in G_m$. But since 
$\BBP\left(\cup_{m\geq 1}G_m\right)=1$,
this proves (\ref{eq:4a}).

Next we turn to the proof of (\ref{eq:4b}).
Since, by the law of large numbers,
$Q_n^*({\cal G}_n)\to 1$, 
as $n\to\infty$, (\ref{eq:4b})
is equivalent to
\be
\liminf_{n\to\infty}\frac{1}{n}\log
        Q^*_n\left(B(X_1^n,D)\cap{\cal G}_n\right)
        \geq -I^*
        \;\;\;\;
        \BBP-
        \mbox{a.s.}
\label{eq:5}
\ee
Choose and fix one of the (almost all)
realizations $x_1^\infty$ of $\BBP$
for which
$$\hat{P}_{x_1^n}(a)\to P(a),
\;\;\;\;\mbox{for all}\; a\in A.$$
Let $\epsilon_1\in(0,\delta)$ arbitrary,
and choose and fix $N$ large enough
so that 
\be
|\hat{P}_{x_1^n}(a)-P(a)|<\epsilon_1P(a)
\;\;\;\;\mbox{for all}\; a\in A,\; n\geq N.
\label{eq:lln}
\ee
Let $a_1,a_2,\ldots, a_m$ denote the elements
of $A$, write $n_0=0$,
$$n_i=n\hat{P}_{x_1^n}(a_i),\;\;i=1,2,\ldots,m$$
and $N_j=\sum_{k=0}^j n_k$, $j=0,1,\ldots,m$.
For $n\geq N$, writing 
$Y_1^n = (Y_1,Y_2,\ldots,Y_n)$
for a vector of random variables with distribution
$Q_n^*$, we have that 
$Q^*_n\left(B(x_1^n,D)\cap{\cal G}_n\right)$
equals
$$Q^*_n\left\{\frac{1}{n}\sum_{i=1}^n\rho(x_i,Y_i)\leq D
	\;\;\mbox{and}\;\;
	\frac{1}{n}\sum_{i=1}^n\IND_{\{Y_i=b\}}
	\leq Q^*(b)+\delta,
	\;\;\;\forall\, b\in A
	\right\}
	\hspace{1.6in}
	$$
$$=Q_n^*\left\{
	\sum_{i=1}^m\frac{n_i}{n}\frac{1}{n_i}
		\sum_{j=N_{i-1}+1}^{N_{i}}
		\rho(a_i,Y_j)\leq D
		\;\;\mbox{and}\;\;
	\sum_{i=1}^m\frac{n_i}{n}\frac{1}{n_i}
                \sum_{j=N_{i-1}+1}^{N_{i}}
		\IND_{\{Y_j=b\}}
		\leq Q^*(b)+\delta,
        	\;\;\;\forall\,b\in A
	\right\}$$
where we have used the fact that the $Y_i$ are IID 
(and hence exchangeable) to rewrite $x_1^n$ as 
consisting of $n_1$ $a_1$'s followed 
$n_2$ $a_2$'s, and so on. Let
$\gamma_i=P(a_i)\sum_{b\in A}W^*(b|a_i)\rho(a_i,b)$
for $i=1,2,\ldots,m$.
Recalling that, by the choice of $W^*$,
$\sum_i\gamma_i=E_{W^*}\rho(X,Y)\leq D,$
and that $Q^*$ is the $Y$-marginal of $W^*$,
the above probability is bounded below by
$$\prod_{i=1}^m
Q_{n_i}^*\left\{
	\frac{n_i}{n}\frac{1}{n_i}
		\sum_{j=1}^{n_{i}}
		\rho(a_i,Y_j)\leq \gamma_i
		\;\;\mbox{and}\;\;
	\frac{n_i}{n}\frac{1}{n_i}
                \sum_{j=1}^{n_{i}}
		\IND_{\{Y_j=b\}}
		\leq P(a_i)[W^*(b|a_i)+\delta],
        	\;\;\;\forall\, b\in A
	\right\}.$$
Writing
$\Gamma_i=\gamma_i/[P(a_1)(1+\epsilon_1)]$,
$i=1,2,\ldots,m$ and 
using (\ref{eq:lln}), 
this is in turn bounded below by
$$\prod_{i=1}^m
Q_{n_i}^*\left\{
	\frac{1}{n_i}
		\sum_{j=1}^{n_{i}}
		\rho(a_i,Y_j)\leq \Gamma_i
		\;\;\mbox{and}\;\;
	\frac{1}{n_i}
                \sum_{j=1}^{n_{i}}
		\IND_{\{Y_j=b\}}
		\leq \frac{W^*(b|a_i)+\delta}{1+\epsilon_1},
        	\;\;\;\forall\, b\in A
	\right\}
	\hspace{0.6in}$$
\vspace*{-0.2in}
\be
	=\;\prod_{i=1}^m
	Q_{n_i}^*\left\{\hat{P}_{Y_1^{n_i}}\in F_i\right\},
	\hspace{3.5in}
	\label{eq:terms}
\ee
where $F_i$ is the collection of
probability mass functions $Q$ on $A$,
$$F_i=F_i(\epsilon_1)=
	\left\{Q\;:\;E_Q[\rho(a_i,Y)]\leq \Gamma_i
	\;\;\mbox{and}\;\;
	Q(b)\leq \frac{W^*(b|a_i)+\delta}{1+\epsilon_1},
                \;\;\;\forall\, b\in A\right\}.$$
We will apply Sanov's theorem
to each one of the terms in (\ref{eq:terms}). Consider 
two cases:
If $\Gamma_i>0$ then $F_i$ is
the closure of its interior (in the
Euclidean topology), so by Sanov's theorem
\be
\liminf_{n_i\to\infty}\frac{1}{n_i}\log
Q_{n_i}^*\left\{\hat{P}_{Y_1^{n_i}}\in F_i\right\}
\geq -\inf_{Q\in F_i} H(Q\|Q^*)
\label{eq:sanov}
\ee
(see Theorem~2.1.10 and Exercise 2.1.16 in 
\cite{dembo-zeitouni:book}).
If $\Gamma_i=0$ then $\gamma_i=0$
and this can only happen if
$W^*(\cdot|a_i)=\IND_{\{a_i\}}(\cdot)$,
in which case $F_i=\{\delta_{a_i}\}$ and
\ben
\frac{1}{n_i}\log 
Q_{n_i}^*\left\{\hat{P}_{Y_1^{n_i}}\in F_i\right\}
\,=\,\log Q^*(a_i)\,=\,-H(\delta_{a_i}\|Q^*)
\een
so (\ref{eq:sanov}) still holds in this case.
Combining the above steps 
(note that each $n_i\to\infty$ as $n\to\infty$),
\ben
\liminf_{n\to\infty} \frac{1}{n}\log
	Q^*_n\left(B(x_1^n,D)\cap{\cal G}_n\right)
&\geq&
	\liminf_{n\to\infty}
	\frac{1}{n}\log \left[ \prod_{i=1}^m
	Q_{n_i}^*\left\{\hat{P}_{Y_1^{n_i}}\in F_i\right\}
	\right]\\
&=& 
	\liminf_{n\to\infty}
	\sum_{i=1}^m \hat{P}_{x_1^{n}}(a_i)
	\frac{1}{n_i}\log
	Q_{n_i}^*\left\{\hat{P}_{Y_1^{n_i}}\in F_i\right\}\\
&\geq&-\sum_{i=1}^m P(a_i)
	\inf_{Q\in F_i} H(Q\|Q^*),
\een
and this holds for $\BBP$-almost any 
$x_1^\infty$. Rewriting the $i$th infimum
above as the infimum over conditional
measures $W(\cdot|a_i)\in F_i$, yields
$$\liminf_{n\to\infty} \frac{1}{n}\log
        Q^*_n\left(B(X_1^n,D)\cap{\cal G}_n\right)
	\geq
	-\inf_{W\in F(\epsilon_1)}
	H(W\|P\!\times\!Q^*)
	\;\;\;\;
	\BBP-\mbox{a.s.}
$$
where $F(\epsilon_1)=
\left\{W\;:\; W_X=P\;\;\mbox{and}\;\;
	W(\cdot|a_i)\in F_i(\epsilon_1),
                \;\;\;\forall\, i=1,2,\ldots,m
	\right\}.$ 
Finally, since $\epsilon_1$ was arbitrary
we can let it decrease to 0 to obtain
\ben
\liminf_{n\to\infty} \frac{1}{n}\log
        Q^*_n\left(B(X_1^n,D)\cap{\cal G}_n\right)
   &\geq&
	\limsup_{\epsilon_1\downarrow 0}
	[-\inf_{W\in F(\epsilon_1)}
        H(W\|P\!\times\!Q^*)]\\
   &\eqa&
	-\inf_{W\in F(0)}H(W\|P\!\times\!Q^*)]\\
   &\geqb&-I^*
	\;\;\;\;
        \BBP-\mbox{a.s.}
\een
This gives (\ref{eq:5}) and completes 
the proof, once we justify steps $(a)$
and $(b)$. Step $(b)$ follows upon
noticing that $W^*\in F(0)$
and recalling that 
$H(W^*\|P\!\times\!Q^*)=I^*$.
Step $(a)$ follows from
the fact that $H(W\|P\!\times\!Q^*)$
is continuous over those $W$
that are absolutely continuous
with respect to $P\!\times\!Q^*$,
and from the observation in Lemma~2
below (verified in the Appendix).
\qed

\vspace{0.1in}

{\sc Lemma 2.}
{\em For all $\epsilon_1>0$ small enough
there exist
$Q_i\in F_i(\epsilon_1)$ such 
that $H(Q_i\|Q^*)<\infty,$
for $1\leq i\leq m$. Therefore,
for all $\epsilon_1>0$ small enough
the exists $W\in F(\epsilon_1)$ 
with $H(W\|P\!\times\!Q^*)<\infty.$}

\vspace{0.1in}

Note that, in the above proof, a somewhat 
stronger result than the one given in Theorem~2 
is established: It is not just demonstrated that 
there exist sets $C_n$ achieving $(i)$ and $(ii)$, 
but that (almost) any sequence of
sets $C_n$ generated by taking approximately 
$e^{nI^*}$ IID samples from 
$\widehat{Q}_n$ will satisfy $(i)$ and $(ii)$.

We also mention that Bucklew 
\cite{bucklew:87} used Sanov's 
theorem to prove the direct part 
of Shannon's data compression 
theorem. The proof of Theorem~2
is similar, except that 
it involves a less direct 
application of Sanov's theorem 
to the sequence of non-product 
measures $\widehat{Q}_n$, and 
the conclusions obtained are
somewhat stronger (pointwise
rather than $L^1$ bounds).
Similarly, in the proof of
Theorem~4, the G\"{a}rtner-Ellis 
theorem from large deviations is 
applied in a manner which 
parallels the approach 
of \cite{bucklew:88}.

\section{The General Case}

Let $A$ be a Polish space (namely, 
a complete, separable metric space)
equipped with its associated Borel 
$\sigma$-algebra ${\cal A}$, and let 
$\BBP$ be a probability measure
on $(A^\NN,{\cal A}^\NN)$.
Also let 
$(\Ahat,\hat{\cal A})$ be a (possibly 
different) Polish space. Given a
nonnegative measurable function
$\rho:A\!\times\!\Ahat\to[0,\infty)$,
define 
$\rho_n:A^n\!\times\!\Ahatn\to[0,\infty)$
as in (\ref{eq:sldistance}). [The reason
for considering $\Ahat$ as possibly 
different from $A$ is motivated by 
the common data compression scenario,
where, in practice, it is often 
the case that original data take values in 
a large alphabet $A$ (for example, Gaussian 
data have $A=\RL$), whereas 
compressed data take values in 
a much smaller alphabet (for example,
Gaussian data on a computer are 
typically quantized to the finite
alphabet $\Ahat$ consisting of all
double precision reals).] 

Let $\{X_n\}$ be a sequence of random
variables distributed according to 
$\BBP$, and for each $n\geq 1$
write $P_n$ for the $n$-dimensional 
marginal distribution of $X_1^n$. 
We say that $\BBP$
is a stationary measure if $X_1^n$ has
the same distribution as $X_{1+k}^{n+k}$,
for any $n,k$.
Let $M:\Ahat\to(0,\infty)$ be a measurable
``mass'' function on $\Ahat$.
To avoid uninteresting technicalities,
we will assume throughout
that $M$ is bounded away from zero,
$M(y)\geq M_*$ for 
some constant $M_*>0$ and all
$y\in\Ahat.$
Next we define the natural analogs 
of the rate functions
$I(P,Q,D)$ and $R(D)$.
For $n\geq 1$, $D\geq 0$ 
and $Q_n$ a probability 
measure on $(\Ahatn,\hat{\cal A}^n)$, let
\be
I_n(P_n,Q_n,D)=\inf_{W_n\in{\cal M}_n(P_n,Q_n,D)}
        H(W_n\|P_n\!\times\!Q_n)
\label{eq:Igendef}
\ee
where $H(\mu\|\nu)$ denotes the relative entropy
between two probability measures
$\mu$ and $\nu$ 
\ben
H(\mu\|\nu)=
        \left\{ \begin{array}{ll}
                        \int d\mu\log\frac{d\mu}{d\nu},&\;\;\;\mbox{when}
                                \frac{d\mu}{d\nu}\;\mbox{exists}\\
                        \infty,                        &\;\;\;\mbox{otherwise}
                \end{array}
        \right.
\een
and where ${\cal M}_n(P_n,Q_n,D)$
consists of all probability measures
$W_n$ on 
$(A^n\!\times\!\Ahatn,
{\cal A}^n\!\times\!\hat{\cal A}^n)$
such that $W_{n,X}$, the first
marginal of $W_n$, is equal to $P_n$, 
the second marginal $W_{n,Y}$ is $Q_n$, 
and $\int \rho_n\, dW_n\leq D$;
if ${\cal M}_n(P_n,Q_n,D)$ is empty,
let $I_n(P_n,Q_n,D)=\infty$.
Then $R_n(D)$ is defined by
\be
R_n(D)=R_n(D;P_n,M) = \inf_{Q_n}
                \left\{ I_n(P_n,Q_n,D) + E_{Q_n}[\log M^n(Y_1^n)]
                \right\},
\label{eq:simpleRgendef}
\ee
where the infimum is over all 
probability measures $Q_n$ on 
$(\Ahatn,\hat{\cal A}^n)$.
Note that since $I_n(P_n,Q_n,D)$
is nonnegative 
and $M$ is bounded away from zero, 
$R_n(D)$ is always well-defined.
Recall also that the mutual information
between two random vectors $X_1^n$ 
and $Y_1^n$ with joint distribution
$W_n$ and corresponding marginals
$P_n$ and Q$_n$, is defined 
by $I(X_1^n;Y_1^n)=H(W_n\|P_n\!\times\!Q_n)$,
so that $R_n(D)$ can alternatively be 
written in a form analogous to (\ref{eq:infoRdef})
in the discrete case:
$$R_n(D)=\inf_{(X_1^n,Y_1^n):\;
X_1^n\sim P_n,\;E\rho_n(X_1^n,Y_1^n)\leq D} \left\{
        I(X_1^n;Y_1^n) + E[\log M^n(Y_1^n)] \right\}.$$
Finally, the rate function $R(D)$ is defined
by
$$R(D)=\lim_{n\to\infty}\frac{1}{n} R_n(D)$$
whenever the limit exists.
Next we state some simple properties of
$R_n(D)$ and $R(D)$, proved in the Appendix.

\vspace{0.1in}

{\sc Lemma 3.}{\em
$(i)$ For each $n\geq 1$, 
$R_n(D)$ is nonincreasing and
convex in $D\geq 0$,
and therefore also continuous
at all $D$ except possibly at
the point
\ben
\Dminn&=&\inf\{D\geq 0\;:\;R_n(D)<+\infty\}.
\een

$(ii)$ If $R(D)$ exists
then it is nonincreasing 
and convex in $D\geq 0$,
and therefore also continuous 
at all $D$ except possibly at 
the point
\ben
\Dmin&=&\inf\{D\geq 0\;:\;R(D)<+\infty\}.
\een

$(iii)$ If $\BBP$ is a stationary measure,
then 
$$R(D)=\lim_{n\to\infty} \frac{1}{n} R_n(D)
=\inf_{n\geq 1}\frac{1}{n} R_n(D)
\;\;\;\mbox{exisits,}$$
and $\Dmin=\inf_n\Dminn$.

$(iv)$ The mutual information $I(X_1^n;Y_1^n)$
is convex in the marginal distribution $P_n$
of $X_1^n$, for a fixed conditional distribution
of $Y_1^n$ given $X_1^n$.
}

\vspace{0.1in}

Next we state analogs of 
Theorems~1 and~2 in the general
case. As before, we are interested
in sets $C_n$ that have large 
blowups but small masses;
since $M$ is bounded away from 
zero we may restrict our attention 
to finite sets $C_n$.

\vspace{0.1in}

{\sc Theorem 3.} {\em Let $C_n\subset \Ahatn$
be an arbitrary finite set
and write 
$D=E_{P_n}[\rho_n(X_1^n, C_n)]$.
Then
\be
\log M^n(C_n)\geq R_n(D).
\label{eq:generalLB}
\ee
If $\BBP$ is a stationary measure,
then for all $n\geq 1$}
\ben
\log M^n(C_n)\geq nR(D).
\een

As will become apparent
from its proof (at the
end of this section), Theorem~3 remains 
true in great generality. The exact same
proof works for arbitrary (non-product)
positive mass functions $M_n$ in place of
$M^n$, and more general distortion measures
$\rho_n$, not necessarily of the form in
(\ref{eq:sldistance}). Moreover,
as long as $R_n(D)$ is well-defined,
the assumption that $M$ is bounded
away from zero is unnecessary.
In that case we can also consider
countably infinite sets $C_n$,
and (\ref{eq:generalLB}) remains
valid as long as $R_n(D)$ is continuous
in $D$
(see Lemma~3).

In the special case when $\BBP$ is
a product measure it is not hard to check
that $R_n(D)=nR(D)$ for all $n\geq 1$,
so we can recover Theorem~1 from Theorem~3.

For Theorem~4 some additional
assumptions are needed. We will assume
that the functions $\rho$ and $\log M$ are
bounded, i.e., that there exist constants 
$\rhomax\geq 0$ and $\Lmax<\infty$
such that $\rho(x,y)\leq \rhomax$ and
$|\log M(y)|\leq \Lmax,$
for all $x\in A$, $y\in\Ahat$.
For $k\geq 1$, we say that 
$\BBP$ is stationary (respectively,
ergodic) in $k$-blocks if the
process
$\{\widetilde{X}_n^{(k)}\;;\;n\geq 0\} 
=\{X_{nk+1}^{(n+1)k}\;;\;n\geq 0\}$ is 
stationary (resp. ergodic). If $\BBP$ 
is stationary then it is stationary in 
$k$-blocks for every $k$. But an ergodic
measure $\BBP$ may not be ergodic in
$k$-blocks. For the second part of the Theorem
we will assume that $\BBP$
is ergodic in blocks, that is, that it
is ergodic in $k$-blocks for all $k\geq 1$.
Also, since $R(D)=\infty$ for $D$ below $\Dmin$,
we restrict our attention to the case 
$D>\Dmin$. Theorem~4 is proved in the 
next section.

\vspace{0.1in}

{\sc Theorem 4.} {\em Assume that 
the functions $\rho$ and $\log M$
are bounded, and that $\BBP$
is a stationary ergodic measure.
For any $D>\Dmin$
and any $\epsilon>0$,
there is a sequence of sets $\{C_n\}$
such that:
\ben
&(i)&\hspace{0.2in}
        \frac{1}{n}\log M^n(C_n)\leq R(D)+\epsilon
        \hspace{0.3in} \mbox{for all $n\geq 1$}\\
&(ii)&\hspace{0.2in}
	P_n([C_n]\subD)\to 1
	\hspace{0.3in}\mbox{as}\;n\to\infty.
\een
If, moreover, $\BBP$ is ergodic 
in blocks, there are sets 
$\{C_n\}$ that satisfy $(i)$ and 
\ben
&(iii)&\hspace{0.2in}
        \rho_n(X_1^n,C_n)\leq D
        \hspace{0.3in}
        \mbox{eventually, }
        \BBP-\mbox{a.s.}
	\;\;\;
\een
}

\vspace*{-0.1in}

{\sc Remark 3.}
A corresponding version of the
asymptotic form of Theorems~1 and~2
given in Remark~1 of the previous
section can also be derived here,
and it holds for every stationary
ergodic $\BBP$. 

\vspace{0.1in}

{\sc Remark 4.}
The assumptions on the boundedness 
of $\rho$ and $\log M$ are made for
the purpose of technical convenience,
and can probably be relaxed to 
appropriate moment conditions.
Similarly, the assumption that
$M^n$ is a product measure can
be relaxed to include sequences
of measures $M_n$ that have rapid 
mixing properties. Finally,
the assumption that $\BBP$
is ergodic in blocks is not as severe
as it may sound. For example,
it is easy to see that any weakly
mixing measure (in the ergodic-theoretic
sense -- see \cite{petersen:book})
is ergodic in blocks.

\vspace{0.1in}


{\em Proof of Theorem~3}:
Given an arbitrary $C_n$,
let $\phi_n:A^n\to C_n$
be defined as in the
proof of Theorem~1.
For $X_1^n\sim P_n$ define
$Y_1^n=\phi_n(X_1^n)$,
write $Q_n$ for the (discrete)
distribution of $Y_1^n$,
and
$W_n(dx_1^n,dy_1^n) =
P_n(dx_1^n)\delta_{\phi_n(x_1^n)}(dy_1^n)$
for the joint distribution
of $(X_1^n,Y_1^n)$.
Then $E_{W_n}[\rho_n(X_1^n,Y_1^n)]=D$,
and by Jensen's inequality applied as
in the discrete case
\ben
\log M^n(C_n)
&\geq& \sum_{y_1^n\in C_n}Q_n(y_1^n)\log
        \frac{M^n(y_1^n)}{Q_n(y_1^n)}\\
&=& \int\,dW_n(x_1^n,y_1^n)\log
        \frac{dW_n(x_1^n,y_1^n)}{d(P_n\!\times\!Q_n)}
        +\sum_{y_1^n\in C_n}Q_n(y_1^n)\log M^n(y_1^n)\\
&=& I(X_1^n;Y_1^n)+E_{Q_n}[\log M^n(Y_1^n)].
\een
By the definition of $R_n(D)$,
this is bounded below by $R_n(D)$.
The second part follows immediately
from the fact that $R_n(D)\geq nR(D)$,
by Lemma~3~$(ii)$.
\qed

\section{Proof of Theorem 4}
The proof of the Theorem is given in 3 steps.
First we assume that $\BBP$ is ergodic
in blocks, and for any $D>\Dminone$ we
construct sets $C_n$ satisfying $(i)$ and
$(iii)$ with $R_1(D)$ in place of $R(D)$.
In the second step
(still assuming $\BBP$ is ergodic
in blocks), for each $D>\Dmin$
we construct sets $C_n$ satisfying
$(i)$ and $(iii)$. In Step~3 we
drop the assumption of the ergodicity 
in blocks, and for any $D>\Dmin$
we construct sets $C_n$ satisfying
$(i)$ and $(ii)$.

\subsection{Step 1:}
Let $\BBP$ and $D>\Dminone$ be
fixed, and let an arbitrary $\epsilon>0$ 
be given.  
By Lemma~3 we can choose 
a $D'\in (\Dmin,D)$ such that
$R_1(D')\leq R_1(D)+\epsilon/8$
and a probability measure $Q^*$ 
on $(\Ahat,\hat{\cal A})$ such that
\be
I^*+L^*\bydef I_1(P_1,Q^*,D')
	+E_{Q^*}[\log M(Y)]
		\leq R_1(D)+\epsilon/8
			\leq R_1(D)+\epsilon/4.
\label{eq:first-order}
\ee
Also we can pick a $W^*\in {\cal M}_1(P_1,Q^*,D')$
such that
\be
H(W^*\|P_1\!\times\!Q^*)\leq I^*+\epsilon/4.
\label{eq:Wchoice}
\ee
As in the proof of Theorem~2,
for $n\geq 1$, write $Q_n^*$ for the product 
measure $(Q^*)^n$, and 
define
$${\cal H}_n=\left\{y_1^n\in \Ahatn\;:\;
	\frac{1}{n}\sum_{i=1}^n
	\log M(y_i)\leq L^*+\epsilon/4 
\right\}.$$
Let $\widetilde{Q}_n$ be the measure $Q_n^*$
conditioned on ${\cal H}_n$,
$\;\widetilde{Q}_n(F)=
Q_n^*(F\cap {\cal H}_n)/Q_n^*({\cal H}_n),\;$
for $\;F\in \hat{\cal A}^n.$
For each $n\geq 1$, let
$\{Y(i)=(Y_1(i),Y_2(i),\ldots,Y_n(i))\;;\;i\geq 1\}$
be IID random vectors
$Y(i)\sim \widetilde{Q}_n$, and define
$$C_n=\{Y(i)\;:\;1\leq i\leq e^{n(I^*+\epsilon/2)}\}.$$
By the definition of ${\cal H}_n$, 
any $y_1^n\in {\cal G}_n$ has
$M^n(y_1^n)\leq e^{n(L^*+\epsilon/4)}$,
so by (\ref{eq:first-order})
$$ M^n(C_n)\leq e^{n(I^*+\epsilon/2)}e^{n(L^*+\epsilon/4)}
        \leq e^{n(R_1(D)+\epsilon)}$$
and $(i)$ of the Theorem is satisfied
with $R_1(D)$ in place of $R(D)$.
Let $X_1^n$ be a random vector
with distribution $P_n$, and, as
in the proof of Theorem~2, let
$i_n$ be the index of the first $Y(i)$
that matches $X_1^n$ within $\rho_n$-distortion $D$.
To verify $(iii)$ we will show that
\ben
i_n\leq e^{n(I^*+\epsilon/2)}
\;\;\;\;\mbox{eventually, }\BBP\!\times\!\RN-
        \mbox{a.s.}
\een
where $\RN=\prod_{n\geq 1}(\widetilde{Q}_n)^\NN$,
and this will follow from 
the following two statements:
\be
\limsup_{n\to\infty}\frac{1}{n}\log\left[i_n\,
        \widetilde{Q}_n(B(X_1^n,D))\right]
        \leq 0\;\;\;\;
        \BBP\!\times\!\RN-
        \mbox{a.s.}
	\label{tmp}\\
\liminf_{n\to\infty}\frac{1}{n}\log
        \widetilde{Q}_n(B(X_1^n,D))
        \geq -(I^*+\epsilon/4)
        \;\;\;\;
        \BBP-\mbox{a.s.}
	\label{eq:asbefore}
\ee
The proof of (\ref{tmp}) is exactly
the same as the proof of (\ref{eq:4a}) in
the proof of Theorem~2.
To prove (\ref{eq:asbefore}),
first note that
by the law of large numbers
$Q_n^*({\cal H}_n)\to 1$,
as $n\to\infty$, so 
(\ref{eq:asbefore})
is equivalent to
\be
\liminf_{n\to\infty}\frac{1}{n}\log
        Q^*_n\left(B(X_1^n,D)\cap{\cal H}_n\right)
        \geq -(I^*+\epsilon/4)
        \;\;\;\;
        \BBP-
        \mbox{a.s.}
\label{eq:new4}
\ee
Let $Y_1,Y_2,\ldots$ be IID random
variables with common distribution $Q^*$.
For any
realization $x_1^\infty$ of $\BBP$,
define the random vectors 
$\xi_i$ and $Z_n$ by
\ben
\xi_i&=&\left(\rho(x_i,Y_i),\,\log M(Y_i)\right),
	\hspace{0.4in}
	i\geq 1\\
Z_n  &=&\frac{1}{n}\sum_{i=1}^n\xi_i,
	\hspace{0.9in}
        n\geq 1.
\een
Also let $\LA_n(\la)$ be the log-moment
generating function of $Z_n$,
$$\LA_n(\la)
	=\log E\left[
	e^{(\las,Z_n)}
	\right],\;\;\;\;\;
	\la=(\la_1,\la_2)\in\RL^2,
	$$
where $(\cdot,\cdot)$ denotes the usual
inner product in $\RL^2$.
Then for $\BBP$-almost any
$x_1^\infty$, by the
ergodic theorem,
\be
\frac{1}{n}\LA_n(n\la)
&=&	\frac{1}{n}\log E\left[
		e^{\sum_{i=1}^n(\las,\xi_i)}
	\right]
	\nonumber\\
&=&	\frac{1}{n}\sum_{i=1}^n\log E_{Q^*}\left[
		e^{\las_1\rho(x_i,Y)+\las_2\log M(Y)}
        \right]
	\nonumber\\
&\to&	E_{P_1}\left\{
		\log E_{Q^*}\left[
		e^{\las_1\rho(X,Y)+\las_2\log M(Y)}
		\right]
	\right\}
	\label{eq:LA1}
\ee
where $X$ and $Y$ above are independent random 
variables with distributions $P_1$ and $Q^*$,
respectively. Next we will need the following 
lemma. Its proof is a simple application of
the dominated convergence theorem, using the
boundedness of $\rho$ and $\log M$.

\vspace{0.1in}

{\sc Lemma 4.} {\em For $k\geq 1$ and probability
measures $\mu$ and $\nu$ on 
$(A^k,{\cal A}^k)$ and
$(\Ahatk,\hat{\cal A}^k)$,
respectively, define
$$\LA_{\mu,\nu}(\la)=\int
	\log 
	\left\{
	\int\left[
	\exp\left(
	\la_1\rho_k(x_1^k,y_1^k)+\la_2\frac{1}{k}\log M^k(y_1^k)
	\right)\right]d\nu(y_1^k)\right\}d\mu(x_1^k),
$$
for $\la=(\la_1,\la_2)\in\RL^2$.
Then $\LA_{\mu,\nu}$ is convex, finite,
and differentiable for all $\la\in\RL^2$.
}

\vspace{0.1in}

From Lemma~4 we have that the limiting
expression in (\ref{eq:LA1}), which 
equals $\LA_{P_1,Q^*}$, is finite
and differentiable everywhere. Therefore
we can apply the 
G\"{a}rtner-Ellis theorem 
(Theorem~2.3.6 in
\cite{dembo-zeitouni:book})
to the sequence of random
vectors $Z_n$, along
$\BBP$-almost any
$x_1^\infty$, to get
\be
\liminf_{n\to\infty}\frac{1}{n}\log
        Q^*_n\left(B(x_1^n,D)\cap{\cal H}_n\right)
=
\liminf_{n\to\infty}\frac{1}{n}\log \PR(Z_n\in F)
\geq 
-\inf_{z\in F}
\LA^*(z)
\;\;\;\;\BBP-\mbox{a.s.}
\label{eq:GE}
\ee
where 
$F=\{z=(z_1,z_2)\in\RL^2\;:\;z_1<D,\;z_2<L^*+\epsilon/4\}$
and
$$\LA^*_{P_1,Q^*}(z)=\sup_{\las\in\RL^2}
	[(\la,z)-\LA_{P_1,Q^*}(\la)]
	$$ 
is the Fenchel-Legendre transform
of $\LA_{P_1,Q^*}(\la)$. 
Recall our
choice of $W^*$ in (\ref{eq:Wchoice}).
Then for any
bounded measurable function
$\phi:\Ahat\to\RL$ and any fixed
$x\in A$,
$$H(W^*(\cdot|x)\|Q^*(\cdot)) 
	\geq \int \phi(y)dW^*(y|x) -\log \int e^{\phi(y)} 
dQ^*(y)$$
(see, e.g., Lemma~6.2.13 in 
\cite{dembo-zeitouni:book}).
Fixing $x\in A$ and $\la\in\RL^2$
for a moment, take $\phi(y)=\la_1\rho(x,y)+\la_2\log M(y)$,
and integrate both sides $dP_1(x)$
to get
$$H(W^*\|P_1\!\times\!Q^*)
\geq\la_1 E_{W^*}(\rho)+\la_2E_{Q^*}[\log M(Y)]
	-\LA_{P_1,Q^*}(\la).$$
Taking the supremum over all $\la\in\RL^2$
and recalling (\ref{eq:Wchoice}) this
becomes
$$I^*+\epsilon/4\geq
H(W^*\|P_1\!\times\!Q^*)
\geq \LA^*_{P_1,Q^*}(D^*,L^*)$$
where
$D^*=\int\rho\,d W^*\leq D'<D$,
so
$$I^*+\epsilon/4\geq
\inf_{z\in F}
\LA^*_{P_1,Q^*}(z).$$
Combining this with the
bound (\ref{eq:GE})
yields
(\ref{eq:new4}) as required,
and completes the proof of this step.

\subsection{Step 2:}
Let $\BBP$ and $D>\Dmin$ be
fixed, and an arbitrary $\epsilon>0$
be given. 
By Lemma~3 we can pick $k\geq 1$ 
large enough so that $\Dmink<D$ and
$(1/k)R_k(D)\leq R(D) +\epsilon/8.$
This step consists of essentially repeating
the argument of Step~1 along blocks of length 
$k$. Choose a $D'\in(\Dmink,D)$ such that
\be
\frac{1}{k}R_k(D')\leq\frac{1}{k}R_k(D)+\epsilon/16,
\label{eq:Dchoice}
\ee
and a probability measure $Q_k^*$
on $(\Ahatk,\hat{\cal A}^k)$ achieving
\be
I_k^*+L_k^*\bydef
\frac{1}{k}I_k(P_k,Q_k^*,D')
+\frac{1}{k} E_{Q_k^*}[\log M^k(Y_1^k)]
\leq 
\frac{1}{k}R_k(D'),
\label{eq:Qkchoice}
\ee
so that
\be
I_k^*+L_k^*\leq R(D)+\epsilon/4.
\label{eq:kth-order}
\ee
Also pick a $W_k^*\in
{\cal M}_k(P_k,Q_k^*,D')$
such that
\be
\frac{1}{k}H(W_k^*\|P_k\!\times\!Q_k^*)\leq I_k^*+\epsilon/4.
\label{eq:Wkchoice}
\ee
For any $n\geq 1$ write $n=mk+r$ for integers
$m\geq 0$ and $0\leq r<k$, and define
$${\cal H}_{n,k}=\left\{y_1^n\in \Ahatn\;:\;
        \frac{1}{n}\sum_{i=1}^n
        \log M(y_i)\leq L_k^*+\epsilon/4
\right\}.$$
Write $Q^*_{n,k}$ for the measure
$$\left[\prod_{i=1}^m Q_k^*\right]
\!\times\!
[Q_k^*]_r,$$
where $[Q_k^*]_r$ denotes the restriction 
of $Q_k^*$ to $(\hat{A}^r,\hat{\cal A}^r)$,
and let $\widetilde{Q}_{n,k}$ be the 
measure $Q_{n,k}^*$
conditioned on ${\cal H}_{n,k}$.
For each $n\geq 1$, let
$\{Y(i)=(Y_1(i),Y_2(i),\ldots,Y_n(i))\;;\;i\geq 1\}$
be IID random vectors
$Y(i)\sim \widetilde{Q}_n$, and let
$C_n$ consist of the first
$e^{n(I_k^*+\epsilon/2)}$ of them.
As before, by the definitions of 
${\cal H}_{n,k}$ and $C_n$,
and using (\ref{eq:kth-order}),
it easily follows that
$$\frac{1}{n}\log M^n(C_n)\leq R(D)+\epsilon$$
so $(i)$ of the Theorem is satisfied.
Let $Y_1,Y_2,\ldots,Y_n$ be 
distributed according to $Q_{n,k}^*$, and note
that the random vectors $Y_{ik+1}^{(i+1)k}$
are IID with distribution $Q_k^*$
(for $i=0,1,\ldots,m-1$). Therefore,
as $n\to\infty$, by the law of large 
numbers we have that with probability 1:
\be
\frac{1}{n}\sum_{i=1}^n\log M(Y_i)
\leq\left(\frac{m}{n}\right)
\frac{1}{m}\sum_{i=0}^{m-1}\log M^k(Y_{ik+1}^{(i+1)k})
	+\frac{k\Lmax}{n}
	\to L_k^*.
\label{eq:k-lln}
\ee
Following the same steps as before,
to verify $(iii)$ it suffices to show that
\ben
\liminf_{n\to\infty}\frac{1}{n}\log
        \widetilde{Q}_{n,k}(B(X_1^n,D))
        \geq -(I_k^*+\epsilon/4)
        \;\;\;\;
        \BBP-\mbox{a.s.}
\een
and, in view of
(\ref{eq:k-lln}), this reduces to
\be
\liminf_{n\to\infty}\frac{1}{n}\log
        Q^*_{n,k}\left(B(X_1^n,D)\cap{\cal H}_{n,k}\right)
        \geq -(I_k^*+\epsilon/4)
        \;\;\;\;
        \BBP-
        \mbox{a.s.}
\label{eq:9}
\ee
For an arbitrary realization $x_1^\infty$ from $\BBP$
and with $Y_1^n$ as above, consider blocks of length $k$.
For $i=0,1,\ldots,m-1,$ we write
$$\widetilde{Y}_i^{(k)}=Y_{ik+1}^{(i+1)k}
\;\;\;\;\mbox{and}\;\;\;\;
\widetilde{x}_i^{(k)}=x_{ik+1}^{(i+1)k}$$
so that the probability
$Q^*_{n,k}\left(B(X_1^n,D)\cap{\cal H}_{n,k}\right)$
can be written as
\ben
Q^*_{n,k}
  \left\{
	\left(\frac{mk}{n}\right)\frac{1}{m}\sum_{i=0}^{m-1}
	\rho_k(\widetilde{Y}_i^{(k)},\widetilde{x}_i^{(k)})
	+\frac{r}{n}\rho_r(Y_{n-r+1}^n,x_{n-r+1}^n)\leq D
	\right.
        \hspace{1.0in}\\
    \;\;\mbox{and}\;\;
	\left.
	\left(\frac{mk}{n}\right)\frac{1}{m}\sum_{i=0}^{m-1}
	\frac{1}{k}\log M^k(\widetilde{Y}_i^{(k)})
	+\frac{1}{n}\log M^r(Y_{n-r+1}^n)\leq L_k^*+\epsilon/4
  \right\}.
\een
Since we assume $\rho(x,y)\leq\rhomax$ 
and $|\log M(y)|\leq\Lmax$ for all $x\in A,\;y\in\Ahat$,
then for all $n$ large enough 
(uniformly in $x_1^\infty$)
the above probability is bounded below by
$$(Q_{k}^*)^m
  \left\{
	\frac{1}{m}\sum_{i=0}^{m-1}
	\rho_k(\widetilde{Y}_i^{(k)},\widetilde{x}_i^{(k)})
	\leq D' +\epsilon/8
     \;\;\mbox{and}\;\;
	\frac{1}{m}\sum_{i=0}^{m-1}
	\frac{1}{k}\log M^k(\widetilde{Y}_i^{(k)})
	\leq L_k^*+\epsilon/8
  \right\}.$$
Now we are in the same situation as in 
the previous step, with the IID random
variables $\widetilde{Y}_i^{(k)}$ in place
of the $Y_i$, the ergodic process
$\{\widetilde{X}_i^{(k)}\}$ in place
of $\{X_i\}$, and $D'+\epsilon/8$ in place 
of $D$. Repeating the same argument as in
Step~1 and invoking Lemma~4 and the 
G\"{a}rtner-Ellis theorem, 
\be
\liminf_{n\to\infty}\frac{1}{n}\log
        Q^*_{n,k}\left(B(X_1^n,D)\cap{\cal H}_{n,k}\right)
        \geq -\inf_{z_1<D'+\epsilon/8,\;z_2<L_k^*+\epsilon/8}
	\LA^*_{k}(z_1,z_2)
        \;\;\;\;
        \BBP-
        \mbox{a.s.}
\label{eq:10}
\ee
where, in the notation of 
Lemma~4, $\LA^*_{k}(z)$ is
the Fenchel-Legendre transform of
$\LA_{P_k,Q_k^*}(\la)$.
Recall our choice of $W_k^*$ in
(\ref{eq:Wkchoice})
and write $D_k^*=\int\rho_k\,dW_k^*\leq D'$.
Then by an application of 
Lemma~6.2.13
from 
\cite{dembo-zeitouni:book}
together with 
(\ref{eq:Wkchoice}) we get that
$$I^*_k+\epsilon/4\geq
\frac{1}{k}H(W_k^*\|P_k\!\times\!Q_k^*)
\geq \LA^*_{k}(D^*,L_k^*),$$
and this together with (\ref{eq:10}) 
proves (\ref{eq:9}), concluding this
step.

\subsection{Step 3:}
In this part we invoke the ergodic
decomposition theorem to remove the
assumption that $\BBP$ is ergodic in
blocks. Although somewhat more 
delicate, the following argument is 
very similar to Berger's proof of the
abstract coding theorem; see pp.~278-281
in \cite{berger:book}.

As in Step~2, let $\BBP$ and $D>\Dmin$ 
be fixed, and let an $\epsilon>0$
be given. Pick $k\geq 1$ large enough 
so that $\Dmink<D$ and 
$\frac{1}{k}R_k(D)\leq R(D) +\epsilon/8,$
and pick $D'\in(\Dmink,D)$ such that
(\ref{eq:Dchoice}) holds. 
Also choose $Q_k^*$ and 
$W_k^*$
as in Step~2 so that
(\ref{eq:Qkchoice}),
(\ref{eq:kth-order})
and (\ref{eq:Wkchoice}) all hold.

Let $\Omega=(A^k)^\NN$, ${\cal F}=({\cal A}^k)^\NN$,
and note that there is a natural 1-1 correspondence
between sets in $F\in{\cal A}^\NN$ and sets in 
$\widetilde{F}\in({\cal A}^k)^\NN$: Writing
$\widetilde{x}_i=x_{ik+1}^{(i+1)k}$,
\be
\widetilde{F}=\{\widetilde{x}_1^\infty\;:\;x_1^\infty\in F\}.
\label{eq:blocking}
\ee
Let $\mu$ be the stationary measure on 
$(\Omega,{\cal F})$ describing
the distribution of the ``blocked'' process
$\{\widetilde{X}_i=X_{ik+1}^{(i+1)k}\;;\;i\geq 0\}$,
where, since $k$ is fixed throughout 
the rest of the proof, we have dropped 
the superscript in $\widetilde{X}_i^{(k)}$.
Although $\mu$ may not be ergodic,
from the ergodic decomposition theorem
we get the following information
(see pp.~278-279 in \cite{berger:book}).

\vspace{0.1in}


{\sc Lemma 5.}{\em There is an integer $k'$ dividing
$k$, and probability measures
$\mu_i,\;i=0,1,\ldots,k'-1$ on $(\Omega,{\cal F})$
with the following properties:

$(i)$ $\;\mu=(1/k')\sum_{i=0}^{k'-1}\mu_i$.

$(ii)$ Each $\mu_i$ is stationary and ergodic.

$(iii)$ For each $i$, let $\BBP^{(i)}$ denote
the measure on $(A^\NN,{\cal A}^\NN)$ induced
by $\mu_i$:
$$\BBP^{(i)}(F)=\mu_i(\widetilde{F}),
\;\;\;\;F\in{\cal A}^\NN$$
[recall the notation
of (\ref{eq:blocking})].
Then $\BBP=(1/k')\sum_{i=0}^{k'-1}\BBP^{(i)},$
and each $\BBP^{(i)}$ is stationary in $k'$-blocks
and ergodic in $k'$-blocks. 

$(iv)$ For each $0\leq i\leq k'$
and $j\geq 0$, the distribution that 
$\BBP^{(i)}$ induces on the process 
$\{X_{j+n}\;;\;n\geq 1\}$ is 
$\BBP^{(i+j\,\mbox{\rm\scriptsize mod}\,k')}$.
}

\vspace{0.1in}

For each $i=0,1,\ldots,k'-1$,
let $\mu_{i,1}$ denote the first-order 
marginal of $\mu_i$ and write
$R(D|i)=R_1(D;\mu_{i,1},\widetilde{M})$ 
for the first-order rate function of
the measure $\mu_i$, with respect to the
distortion measure $\rho_k$, and with
mass function $\widetilde{M}=M^k$.
Since $W_k^*$ chosen as above has its
$A^k$-marginal equal to $P_k$ we can
write it as $W_k^*=V^*_k\circ P_k$ where
$V^*_k(\cdot|X_1^n)$ denote the regular 
conditional probability distributions.
Write $P_k^{(i)}$ for the $k$-dimensional 
marginals of the measures $\BBP^{(i)}$, 
and define probability measures $W_k^{(i)}$ 
on $(A^n\!\times\Ahatn,{\cal A}^n\!\times
\!\hat{\cal A}^n)$
by
$W_k^{(i)}=V^*_k\circ P^{(i)}_k.$
Let $D_i=\int\rho_k\,dW_k^{(i)}$
so that by Lemma~5~$(iii)$,
\be
\frac{1}{k'}\sum_{i=0}^{k'-1}
	D_i
	=\int\rho_k\,dW_k^*\leq D'.
\label{eq:3-4}
\ee
Similarly, writing
$Q_k^{(i)}$ for the $\hat{A}^k$-marginal
of $W_k^{(i)}$ and applying Lemma~5~$(iii)$,
\be
\frac{1}{k'}\sum_{i=0}^{k'-1}
	\int\log M^k(y_1^k)\,dQ_k^{(i)}(y_1^k)
	=\int\log M^k(y_1^k)\,dQ_k^*(y_1^k)
\label{eq:3-6}
\ee
and using the convexity
of mutual information from 
Lemma~3~$(iv)$,
\be
\frac{1}{k'}\sum_{i=0}^{k'-1}
        H(W_k^{(i)}\|P_k^{(i)}\!\times\!Q_k^{(i)})
        \leq H(W_k^*\|P_k\!\times\!Q_k^*).
\label{eq:3-5}
\ee
For $N\geq 1$ large enough we can use
result of Step~1 to get $N$-dimensional
sets $B_i$ that almost-cover $(\hat{A}^k)^N$
with respect to $\mu_i$. Specifically,
consider $N$ large enough so that
\be
\frac{\max\{\rhomax,\;\Lmax,\;1\}}{kN}
< \min\{\epsilon/8,\;(D-D')/2\}.
    \label{eq:3-6b}
\ee
For any such $N$, by the result of Step~1
we can choose sets $B_i\subset(\hat{A}^k)^N$
such that, for each $i$,
\be
\mu_i\left([B_i]\subDi\right)
	&\geq& 1-\epsilon_N,
	\;\;\;\;\mbox{where}\;\epsilon_N\to 0
	\;\mbox{as}\;N\to\infty,\;\;\mbox{and}
	\label{eq:3-7}\\
\widetilde{M}^N(B_i)
	&\leq&\exp\{N(R(D_i|i)+\epsilon/8)\}.
	\label{eq:3-8}
\ee
Now choose and fix an arbitrary $y^*\in{\hat A}$,
and for $n=k'(Nk+1)$ define new sets 
$B_i^*\subset \hat{A}^n$ by
\ben
B_i^*=\prod_{j=0}^{k'-1}
	\left[
	B_{i+j\,\mbox{\rm\scriptsize mod}\,k'}
	\!\times\!\{y^*\}
	\right],
\een
where $\prod$ denotes the cartesian product.
Then, by ({\ref{eq:3-6b}), for any $x_1^n$,
$$\rho_n(x_1^n,B_i^*)<\frac{D-D'}{2}
	+\frac{1}{k'}\sum_{j=0}^{k'-1}
	\rho_{kN}
 	  \left(
		x_{j(kN+1)+1}^{j(kN+1)+kN},
		B_{i+j\,\mbox{\rm\scriptsize mod}\,k'}
	  \right),$$
so by a simple union bound,
\be
\BBP^{(i)}\left([B^*_i]\subD\right)
&\geqa&
	1-\sum_{j=0}^{k'-1}
	\left[ 1-
	  \BBP^{(i+j\,\mbox{\rm\scriptsize mod}\,k')}
	  \left(
	    [B_{i+j\,\mbox{\rm\scriptsize mod}\,k'}]\subD
	  \right)
	\right]
	\nonumber\\
&\eqb&
	1-\sum_{i=0}^{k'-1}
	\left[ 1-
	  \mu_i\left([B_i]\subDi\right)\
	\right]
	\nonumber\\
&\geqc&
1-k'\epsilon_N,
	\label{eq:check1}
\ee
where we used (\ref{eq:3-4}) in $(a)$,
Lemma~5~$(iv)$ in $(b)$,
and (\ref{eq:3-7}) in $(c)$.
Also, using the definition
of $B_i^*$ and the bounds (\ref{eq:3-6b})
and (\ref{eq:3-8}),
\ben
\frac{1}{n}\log M^n(B_i^*)
&\leq&
	\frac{\log M(Y^*)}{kN+1} 
	+ \frac{1}{k'} \sum_{j=0}^{k'-1}
	\left[
	  \frac{1}{kN}\log\widetilde{M}^N
	    (B_{i+j\,\mbox{\rm\scriptsize mod}\,k'})
	\right]\\
&\leq&
	\epsilon/8
	+ \frac{1}{k'} \sum_{j=0}^{k'-1}
        \left[
          \frac{1}{k}(R(D_j|j)+\epsilon/8)
	\right],
\een
but from the definition of $R(D|j)$
and (\ref{eq:3-5}) and (\ref{eq:3-6})
this is 
\be
\frac{1}{n}\log M^n(B_i^*)
&\leq&
	\epsilon/4
	+ \frac{1}{k'} \sum_{j=0}^{k'-1}
	\left[
	  \frac{1}{k}H(W_k^{(j)}\|P_k^{(j)}\!\times\!Q_k^{(j)})
	  +\frac{1}{k}
	  \int\log M^k(y_1^k)\,dQ_k^{(j)}(y_1^k)
	\right]
	\nonumber\\
&\leq&
	I_k^*+L_k^*+\epsilon/2
	\nonumber\\
&\leq&
	R(D)+3\epsilon/4,
	\label{eq:check2}
\ee
where the last inequality follows from 
(\ref{eq:kth-order}).
So in (\ref{eq:check1}) and 
(\ref{eq:check2}) we have shown that,
for {\em all} $i=0,1,\ldots,k'-1$,
\be
\BBP^{(i)}\left([B^*_i]\subD\right)
&\geq&1-k'\epsilon_N\;\;\;\;\mbox{and}\\
\frac{1}{n}\log M^n(B_i^*)
&\leq& R(D)+3\epsilon/4.
\ee
Finally we define sets $C_n\subset\Ahatn$
by
$$C_n=\cup_{i=0}^{k'-1} B_i^*.$$
From the last two bounds above
and (\ref{eq:3-6b}),
the sets $C_n$ have
$$
\frac{1}{n}\log M^n(C_n)
\leq\frac{\log k'}{n}+R(D)+3\epsilon/4
\leq R(D) +\epsilon,$$
and by Lemma~5~$(iii)$,
$$
P_n\left([C_n]\subD\right)
=\frac{1}{k'} \sum_{i=0}^{k'-1}\BBP^{(i)}\left([C_n]\subD\right)
\geq
\frac{1}{k'} \sum_{i=0}^{k'-1}\BBP^{(i)}\left([B^*_i]\subD\right)
\geq 1-\epsilon'_n
$$
where $\epsilon'_n=k'\epsilon_N$ when
$n=k'(Nk+1)$.

In short, we have shown that for any $D>\Dmin$
and any $\epsilon>0$, there exist (fixed)
integers $k$, $k'$ and $N_0$ such that:
\ben
 (+)\;\cases{
         \hspace{0.05in}
         \mbox{There is a sequence of sets $C_n$,
	      for $n=k'(Nk+1)$, $N\geq N_0$, 
              satisfying:
             }&\cr
         \hspace{0.4in}
		(1/n)\log M^n(C_n)\leq R(D)+\epsilon
        	\;\;\mbox{for all $n,\;$ and}
	      &\cr
         \hspace{0.4in}
		P_n\left([C_n]\subD\right)\to 1
        	\;\;\mbox{as}\; n\to\infty. 
	 \cr
        }
\een
Since this is essentially an asymptotic 
result, the restrictions that $N\geq N_0$ 
and $n$ be of the form $n=k'(Nk+1)$ are 
inessential. Therefore they can be easily 
dropped to give $(+)$ for all $n\geq 1$,
that is, to produce a sequence of sets 
$\{C_n\;;\;n\geq 1\}$ satisfying $(i)$ 
and $(ii)$ of Theorem~4.  
\qed

\newpage 
\section*{Appendix}

\hspace{0.2in}
{\em Proof of Remark 1:} In view of part $(i)$ 
of Theorem~2 and the remark that
$P^n\left([C_n]\subD\right)\to 1$,
for every $m\geq 1$ we can choose
a sequence of sets $\{C_n^{(m)}\;;\;n\geq 1\}$ such
that
\ben
\frac{1}{n}\log M^n(C^{(m)}_n)
&\leq& R(D) + \frac{1}{m},
        \hspace{0.3in}
        \mbox{for all }
        m,n\geq 1, \;\mbox{and}\\
P^n\left([C^{(m)}_n]\subD\right)
&\geq& 1-\frac{1}{m},
        \hspace{0.3in}
        \mbox{for all }
        m\geq 1,\; n\geq N(m),
\een
where $N(m)$ is some fixed sequence of integers,
strictly increasing to $\infty$ as $m\to\infty$.
So for each $n\geq 1$ there is a unique $m=m(n)$
such that $N(m)\leq n<N(m+1)$. Since 
$\{N(m)\}$ is strictly increasing, the sequence
$\{m(n)\}$ is nondecreasing and $m(n)\to\infty$
as $n\to\infty$.
Define $C_n^*=C_n^{m(n)}$ for all $n\geq 1$.
From the last two bounds,
\ben
\frac{1}{n}\log M^n(C^*_n)
&\leq& R(D) + \frac{1}{m(n)},
        \hspace{0.3in}
        \mbox{for all }
        n\geq 1, \;\mbox{and}\\
P^n\left([C^{*}_n]\subD\right)
&\geq& \frac{1}{m(n)},
        \hspace{0.3in}
        \mbox{for all }
        n\geq N(m(n)).
\een
But since $n$ is always $n\geq N(m(n))$ by definition,
and $m(n)\to\infty$ as $n\to\infty$, this proves $(i')$
and $(ii')$. Also, since $\rho$ is bounded, 
$(iii')$ immediately follows from $(ii')$.
\qed

\vspace{0.1in}

{\em Proof outline of Lemma~1}:
For part $(i)$
it suffices to consider the case
$I(P,Q,D)<\infty$, so we may assume 
that the set ${\cal M}(P,Q,D)$ is 
nonempty. Since the marginals 
of any $W\in{\cal M}(P,Q,D)$ are
$P$ and $Q$, $W$ is absolutely
continuous with respect to $P\!\times\!Q$,
so $H(W\|P\!\times\!Q)$ is
continuous over $W\in{\cal M}(P,Q,D)$.
Since the sets ${\cal M}(P,Q,D)$ are
compact (in the Euclidean
topology), the infimum in 
(\ref{eq:Idef}) must be achieved.
A similar argument works for
$R(D)$: Combining the two infima 
in its definition, 
\be
R(D) = \inf_{W\in{\cal M}(P,D)}
                \left\{ H(W\|W_{X}\!\times\!W_{Y}) 
			+ E_{W_{Y}}[\log M(Y)]
                \right\},
\label{eq:betterRdef}
\ee
where ${\cal M}(P,D)=\cup_Q {\cal M}(P,Q,D)$.
Since the sets ${\cal M}(P,D)$ are compact, 
the infimum in (\ref{eq:betterRdef})
is achieved by some $W^*\in{\cal M}(P,D)$,
and $Q^*=W^*_Y$ achieves the infimum in
(\ref{eq:simpleRdef}).

For part $(ii)$ recall the assumption 
that for all $a\in A$ there is $b=b(a)$ such
that $\rho(a,b)=0$. If we let
$W(a,b)=P(a)\IND_{\{b=b(a)\}}$,
then $W\in {\cal M}(P,D)$ for any
$D\geq 0$ and from (\ref{eq:betterRdef}),
$R(D)\leq E_{W_Y}[\log M(Y)]<\infty$
for all $D\geq 0$. 
Since the sets ${\cal M}(P,D)$ 
are increasing in $D$,
$R(D)$ is nonincreasing. 
To see that it is convex, 
let $W\in{\cal M}(P,D_1)$ 
and $W'\in{\cal M}(P,D_2)$
arbitrary.
Given $\lambda\in[0,1]$ let $\la'=1-\la$,
and write $V=\la W +\la' W'$.
Then $V\in{\cal M}(P,\la D_1 + \la' D_2)$
and the $Y$-marginal of $V$, 
$V_Y$, is $\la W_Y+\la' W'_Y$.
Recalling (\ref{eq:betterRdef})
and that relative 
entropy is jointly convex in 
its two arguments,
\ben
\lefteqn{R(\la D_1 +\la' D_2)}\\
& & \leq H(V\|V_X\!\times\!V_Y) 
	+ E_{V_Y}[\log M(Y)]\\
& & \leq \la \left\{H(W\|W_X\!\times\!W_Y) 
	+ E_{W_Y}[\log M(Y)]\right\}
	+\la' \left\{
	H(W'\|W'_X\!\times\!W'_Y)
	+ E_{W'_Y}[\log M(Y)]\right\}.
\een
Taking the infimum over all $W\in{\cal M}(P,D_1)$,
$W'\in{\cal M}(P,D_2)$,
and using (\ref{eq:betterRdef})
shows that
$R(D)$ is convex, and since
it is finite for all $D\geq 0$ it
is also continuous.

The proof of $(iii)$ is essentially
identical to that of $(ii)$, using the
definition (\ref{eq:Idef}) in place of
(\ref{eq:betterRdef}). The only difference
is that $I(P,Q,D)$ can 
be infinite, so its convexity 
(and the fact that it is nonincreasing)
imply that it is continuous for $D\geq 0$
except possibly at
$D=\inf\{D\geq 0\;:\;I(P,Q,D)<\infty\}$.

Part $(iv)$ is a well-known
information theoretic fact; see, e.g.,
Lemma~9.4.2 in \cite{gray-90:book}.


For part $(v)$ let $W^*$ achieve the
infimum in (\ref{eq:betterRdef}).
Since relative entropy is nonnegative
we always have
$R(D)\geq \Rmin$,
with equality if and only if
$W^*_Y$ is supported on the set
$A'=\{y\in A\;:\;\log M(y)=\Rmin\}$
and $W^*=W^*_X\!\times\!W^*_Y$. Clearly,
these two conditions are
satisfied if and only if 
$$D\geq \inf \{E_{P\!\times\!Q}[\rho(X,Y)]\;:\;
	Q\;\mbox{supported on}\;A'\},$$
but the right hand side above is exactly equal
to $\Dmax$.
\qed

\vspace{0.1in}

{\em Proof of Lemma~2}:
If $\gamma_i=0$ then,
as discussed in the proof of Theorem~2,
$F_i(\epsilon_1)=\{\delta_{a_i}\}$ for
all $\epsilon_1$
and 
$$H(\delta_{a_i}\|Q^*)=-\log Q^*(a_i)
\leq -\log P(a_i)<\infty.$$
If $\gamma_i>0$ then there must 
exist a $b^*\in A$, $b^*\neq a_i$,
such that $W^*(b^*|a_i)>0$.
Write $\dmax$ for the maximum
of $\sum_b W^*(b|a_j)\rho(a_j,b)$
over all $j=1,\ldots,m$, and
let $\rhomin=\min\{\rho(a,b)\;:\;a\neq b\}$.
For $\alpha\in(0,1)$, let
\ben
Q_i(b)\;=\;
        \cases{
        W^*(a_1|a_i)+\alpha& if $b=a_i$\cr
        W^*(b^*|a_i)-\alpha& if $b=b^*$\cr
        W^*(b|a_i)& otherwise.\cr
              }
\een
Then, for $\epsilon_1$  small enough to
make
$(\delta-\epsilon_1)\rhomin>\epsilon_1\dmax(1+\epsilon_1),$
it is an elementary calculation 
to verify that 
$Q_i\in F_i(\epsilon_1)$ 
and $H(Q_i\|Q^*)<\infty$,
as long as $\alpha$ satisfies
the following conditions:
\ben
\alpha &<& 1-W^*(a_i|a_i)\\
\alpha &<& W^*(b^*|a_i)\\
\frac{\epsilon_1\dmax}{\rhomin}\;<\;\alpha &<& 
	\frac{\delta-\epsilon_1}{1+\epsilon_1}.
\een
Taking $W(a_i,b)=Q_i(b)P(a_i)$ we 
also have $W\in F(\epsilon_1)$.
\qed

\vspace{0.1in}

{\em Proof of Lemma~3}: Since
the sets ${\cal M}_n(P_n,Q_n,D)$ are increasing
in $D$, $R_n(D)$ is nonincreasing in $D$.
Next we claim that relative entropy is 
jointly convex in its two arguments. Let
$\mu$, $\nu$ be two probability measures
over a Polish space $(S,{\cal S})$. In the case when
$\mu$ and $\nu$ both consist of only a
finite number of atoms, the joint
convexity of $H(\mu\|\nu)$ is well-known
(see, e.g., Theorem~2.7.2 in 
\cite{cover:book}). In
general, $H(\mu\|\nu)$ can be written as 
$$H(\mu\|\nu)=\sup_{\{E_i\}}\sum_i\mu(E_i)
\log\frac{\mu(E_i)}{\nu(E_i)}$$
where the supremum is over all finite measurable
partitions of $S$
(see Theorem~2.4.1 in \cite{pinsker:book}).
Therefore $H(\mu\|\nu)$ is the pointwise supremum
of convex functions, hence itself convex.
As in (\ref{eq:betterRdef}),
combining the two infima,
$R_n(D)$ can equivalently 
be written as
\be
R_n(D) = \inf_{W_{n}\in{\cal M}_n(P_n,D)}
                \left\{ H(W_n\|W_{n,X}\!\times\!W_{n,Y}) 
			+ E_{W_{n,Y}}[\log M^n(Y_1^n)]
                \right\}
\label{eq:betterRgen}
\ee
where ${\cal M}_n(P_n,D)=\cup_{Q_n} {\cal M}_n(P_n,Q_n,D)$.
Using this together with the joint convexity
of relative entropy as in the proof of Lemma~1~$(ii)$
shows that $R_n(D)$ is convex. Since it is also
nonincreasing and bounded away from $-\infty$,
$R_n(D)$ is also continuous
at all $D$ except possibly at
the point
\ben
\Dminn&=&\inf\{D\geq 0\;:\;R_n(D)<+\infty\}.
\een
This proves $(i)$. For $(ii)$ notice that
if $R(D)$ exists then it must also be 
nonincreasing and convex in $D\geq 0$
since $R_n(D)$ is; therefore,
it must also be continuous
except possibly at $\Dmin$.

For part $(iii)$, let $m,n\geq 1$ arbitrary,
and let $W_m\in {\cal M}_m(P_m,D)$
and $W_n\in {\cal M}_n(P_n,D)$.
Define a probability measure $W_{m+n}$ on
$(A^n\!\times\!\Ahatn,
{\cal A}^n\!\times\!\hat{\cal A}^n)$ by
$$W_{m+n}(dx_1^{m+n}, dy_1^{m+n})
=W_{m}(dy_1^{m}|x_1^{m})
W_{n}(dy_{m+1}^{m+n}|x_{m+1}^{m+n})
P(dx_1^{m+n}).$$
Notice that $W_{m+n}\in{\cal M}_{m+n}(P_{m+n},D),$
and that, if $(X_1^{m+n},Y_1^{m+n})$ are random
vectors distributed according to $W_{m+n}$, then
$Y_1^m$ and $Y_{m+1}^{m+n}$ are conditionally
independent given $X_1^{m+n}$.
Therefore, 
\ben
R_{m+n}(D)
&\leqa&
	H(W_{m+n}\|W_{m+n,X}\!\times\!W_{m+n,Y}) 
	+ E_{W_{m+n,Y}}[\log M^{m+n}(Y_1^{m+n})]\\
&=& I(X_1^{m+n};Y_1^{m+n})
	+ E_{W_{m+n,Y}}[\log M^{m+n}(Y_1^{m+n})]\\
&\leqb&
	I(X_1^m;Y_1^m)
	+ I(X_{m+1}^{m+n};Y_{m+1}^{m+n})
	+ E_{W_{m,Y}}[\log M^m(Y_1^m)]
	+ E_{W_{n,Y}}[\log M^n(Y_1^n)]
\een
where $(a)$ follows
from (\ref{eq:betterRgen}) and 
$(b)$ follows from the conditional
independence of 
$Y_1^m$ and $Y_{m+1}^{m+n}$ 
given $X_1^{m+n}$
(see, e.g., Lemma~9.4.2 in \cite{gray-90:book}).
So we have shown that $R_{m+n}(D)$ is
bounded above by
$$
	H(W_m\|W_{m,X}\!\times\!W_{m,Y}) 
	+ E_{W_{m,Y}}[\log M^{m}(Y_1^{m})]
	+ H(W_n\|W_{n,X}\!\times\!W_{n,Y}) 
	+ E_{W_{n,Y}}[\log M^{n}(Y_1^{n})],
$$
and taking the infimum over all 
$W_m\in {\cal M}_m(P_m,D)$
and $W_n\in {\cal M}_n(P_n,D)$ yields
\be
R_{m+n}(D)\leq R_{m}(D)+ R_{n}(D).
\label{eq:tmp}
\ee
[Note that in the above argument we 
implicitly assumed that we could find
$W_m\in {\cal M}_m(P_m,D)$
and $W_n\in {\cal M}_n(P_n,D)$; if this
was not the case,
then either
$R_m(D)$ or $R_n(D)$ would be equal to $+\infty$,
and (\ref{eq:tmp}) would still trivially hold.]
Therefore the sequence $\{R_n(D)\}$
is subadditive and it follows that
$\lim_n (1/n)R_n(D)=\inf_n (1/n)R_n(D)$.
From this it is immediate that
$\Dmin=\inf_n\Dminn$.

Part $(iv)$ is a well-known information
theoretic fact; see, e.g., Problem~7.4 in
\cite{berger:book}.
\qed


\end{document}